\documentclass{aiaa-tc}

 \usepackage{overcite}
 \usepackage{graphicx}
 \usepackage{lastpage}
 \usepackage{array}
 \usepackage{fancyhdr}
 \usepackage{setspace}
 \usepackage{ifthen}
 \usepackage{varioref}
 \usepackage{wrapfig}
 \usepackage{threeparttable}
 \usepackage{dcolumn}
 \newcolumntype{d}{D{.}{.}{-1}}
 \usepackage[noprefix]{nomencl}
 \makenomenclature
 \usepackage{subfigure}
 \usepackage{subfigmat}
 \usepackage{fancyvrb}
 \fvset{fontsize=\footnotesize,xleftmargin=2em}
 \usepackage{url}

 \usepackage{color,fancybox}
 \usepackage{colortbl}
 \usepackage{amsmath} 
 \usepackage{amssymb}  
 \usepackage{bm} 
 \usepackage{amsfonts} 
 \usepackage{multirow}

\newcommand{\B}[1]{\ensuremath{\mathcal{B}_{#1}}}
\newcommand{\norm}[1]{\ensuremath{\left\| #1 \right\|}}
\newcommand{\bracket}[1]{\ensuremath{\left[ #1 \right]}}

\newcommand{\parenth}[1]{\ensuremath{\left( #1 \right)}}
\newcommand{\reftab}[1]{\ref{tab:#1}}
\newcommand{\refeqn}[1]{(\ref{eqn:#1})}
\newcommand{\reffig}[1]{\ref{fig:#1}}
\newcommand{\tr}[1]{\mbox{tr}\ensuremath{\negthickspace\bracket{#1}}}
\newcommand{\deriv}[2]{\ensuremath{\frac{\partial #1}{\partial #2}}}

\newcommand{\SO}{\ensuremath{\mathrm{SO(3)}}}
\newcommand{\so}{\ensuremath{\mathfrak{so}(3)}}
\newcommand{\SE}{\ensuremath{\mathrm{SE(3)}}}

\renewcommand{\Re}{\ensuremath{\mathbb{R}}}

\title{Polyhedral Potential and Variational Integrator Computation of the Full Two Body Problem}

\author{
Eugene G. Fahnestock\thanks{Graduate Student, Department of
Aerospace Engineering, University of Michigan, 2008 FXB Building,
1320 Beal Avenue, Ann Arbor, MI 48109, AIAA Student Member.},
Taeyoung Lee\thanksibid{1}, Melvin Leok\thanks{Assistant Professor,
Department of Mathematics, University of Michigan.},
\\
N. Harris McClamroch\thanks{Professor, Department of Aerospace
Engineering, University of Michigan, AIAA Senior Member.}, and
Daniel J. Scheeres\thanks{Associate Professor, Department of
Aerospace Engineering, University of Michigan, AIAA Associate
Fellow.}
\\
{\normalsize\itshape University of Michigan, Ann Arbor, Michigan,
USA} }

 \AIAApapernumber{2006-6289}
 \AIAAconference{AIAA/AAS Astrodynamics Specialist Conference, August 21-24, 2006, and Keystone, Colorado}
 \AIAAcopyright{\AIAAcopyrightB{2006}}

\begin{document}

\maketitle

\begin{abstract}
We present a combination of tools which allows for investigation of
the coupled orbital and rotational dynamics of two rigid bodies with
nearly arbitrary shape and mass distribution, under the influence of
their mutual gravitational potential. Methods for calculating that
mutual potential and resulting forces and moments for a polyhedral
body representation are simple and efficient. Discrete equations of
motion, referred to as the Lie Group Variational Integrator (LGVI),
preserve the structure of the configuration space, \SE, as well as
the geometric features represented by the total energy and the total
angular momentum. The synthesis of these approaches allows us to
simulate the full two body problem accurately and efficiently.
Simulation results are given for two octahedral rigid bodies for
comparison with other integration methods and to show the qualities
of the results thus obtained. A significant improvement is seen over
other integration methods while correctly capturing the interesting
effects of strong orbit and attitude dynamics coupling, in multiple
scenarios.
\end{abstract}

\printnomenclature 

\section{Introduction}

The full two body problem studies the dynamics of two irregular rigid bodies interacting under a mutual potential. The full two body problem arises in numerous engineering and scientific fields. Our focus is on the dynamics of two rigid bodies in space due to their mutual gravitational potential. This depends on both the relative position and the relative attitude of the bodies. Therefore, the translational orbit dynamics and the rotational attitude dynamics are coupled in the full two body problem. For example, the trajectory of a very large spacecraft around the Earth is affected by the attitude of the spacecraft, and the dynamics of a binary asteroid pair are characterized by the non-spherical mass distributions of the two bodies. Recently, interest in the full (two) body problem as applied to binary objects in space has increased, due to evidence that such binary systems are common among the overall asteroid population and form an especially large percentage of near-earth objects, with some estimates as high as 16 percent \cite{Bottke-Icarus-1996,Bottke-Nature-1996,Margot-Science-2002,Merline-ASTIII-2002}. In addition, advances in radar and optical observation methods \cite{Ostro-ASTIII-2002} have allowed for better modelling of the shapes and characteristics of bodies in binary systems, encouraging further study of their dynamics.

Some previous work of relevance includes Maciejewski's presentation
of equations of motion of the full two body problem in inertial and
relative coordinates~\cite{Maciejewski-CMDA-1995}. He also discussed
the existence of relative equilibria. Scheeres derived a stability
condition for the full two body problem~\cite{Scheeres-CMDA-2002},
and he studied the planar stability of an ellipsoid-sphere
model~\cite{Scheeres-NTAC-2003}. Spacecraft motion about binary
asteroids has been discussed using the restricted three body
model~\cite{Scheeres-AASastrospec-2003,Gabern-ICDSDE-2004}, and the
four body model~\cite{Scheeres-DS-2005}.

The mutual gravitational potential of two rigid celestial bodies has
been expressed using spherical
harmonics~\cite{Borderies-celesmech-1978,Braun-thesis-1991}. But,
the harmonic expansion is not guaranteed to converge. Convergence is
shown to be an unstable property of such spherical harmonic
series~\cite{Moritz-book-1980}. By this we mean that an arbitrarily
small change to the mass distribution may cause a previously
convergent series to diverge. Another commonly used approach for
evaluating the mutual gravitational potential is to fill each rigid
body's volume with a distribution of point masses, fixed with
respect to one another, the sum of which equals the respective
body's total mass~\cite{Geissler-Icarus-1996,Ashenberg-JGCD-2005}.
Although the mutual potential obtained for two rigid bodies using
this approach converges to the true gravity field in the limit as
the number of point masses becomes arbitrarily large, there are
significant errors in the computation of gravitational forces from
that mutual potential~\cite{Werner-CMDA-1997}.

These problems are avoided with the approach presented in Ref.
\citen{Werner-CMDA-2005} for calculating mutual gravitational
potential based on polyhedral body models. This method converges at
all points exterior to the two bodies. By representing the bodies as
polyhedra, the flexibility of the formalism over more specialized
representations such as spheres and ellipsoids is increased. A
polyhedral rigid body model can directly include important features
such as craters, caves, or deep clefts where contact binaries meet.
The entire body does not have to be modelled uniformly at a high
resolution. The errors in the mutual potential computation can be
reduced to the level of error in each body's shape determination,
and to the level of discretization chosen for that shape.
Derivatives of this polyhedral mutual potential formulation are
given in Ref. \citen{Fahnestock-AASastrospec-2005} and, with some
corrections and improvement of notation, in Ref.
\citen{Fahnestock-CMDA-2006}. These derivatives determine the forces
and torques exerted by the bodies on each other, for use in either
inertial or relative equations of motion that describe the full body
dynamics.

These full body dynamics arise from Lagrangian and Hamiltonian
mechanics; they are characterized by symplectic, momentum and energy
preserving properties. These geometric features determine the
qualitative behavior of the full body dynamics, and they can serve
as a basis for further theoretical study of the full body problem.
The configuration space of the full body dynamics have a Lie group
structure referred to as the Euclidean group, $\SE$. However,
general numerical integration methods, including the widely used
Runge-Kutta schemes, neither preserve the Lie group structure nor
these geometric properties~\cite{Hairer-book-2000}.

The variational approach~\cite{Marsden-actanum-2001} and Lie group
methods~\cite{Iserles-actanum-2000} provide systematic methods of
constructing structure preserving numerical integrators. The idea of
the variational approach is to discretize Hamilton's principle
rather than the continuous equations of
motion~\cite{Marsden-actanum-2001}. The numerical integrator
obtained from the discrete Hamilton's principle exhibits excellent
energy properties, conserves first integrals, and preserves the
symplectic structure. Lie group methods consist of numerical
integrators that preserve the geometry of the configuration space by
automatically remaining on the Lie
group~\cite{Iserles-actanum-2000}. A Lie group method is explicitly
adopted for the variational integrator in Ref.
\citen{Lee-IEEECCA-2005} and \citen{Lee-CMAME-2005}. This unified
integrator, hereafter referred to as the Lie Group Variational
Integrator (or LGVI for short), is symplectic and momentum
preserving, and it exhibits good total energy behavior for
exponentially long time periods. It also preserves the Euclidian Lie
group structure without the use of local charts, reprojection, or
constraints.

Numerical simulation of the full body problem involves two major
problems; a large computational burden in computing mutual
gravitational forces and moments, and inaccuracy of numerical
integrators. The forces and moments must constantly be reevaluated
with any position change or any orientation change, not only at each
time step but at each sub-step involved in the differencing scheme
behind any general numerical integrator. Therefore such general
numerical integrators amplify the computation cost for finding the
forces and moments. The accuracy of such general integrators also
rapidly degrades as the simulation time increases. They fail to
preserve the conserved quantities such as total energy and angular
momentum, which determine the qualitative behavior of the full body
dynamics. Attitude errors tend to accumulate as a consequence of
numerical errors, and this attitude degradation causes significant
errors in the gravitational force and moment computation.

The unified treatment given in this paper of the polyhedral mutual
potential formulation combined with the LGVI presents a solution to
these two major dynamic simulation problems. Using polyhedral
models, we can approximate irregular bodies to a specified accuracy,
and we can control the computational burden to compute the mutual
gravitational forces and moments by choosing the level of body
discretization and the number of series terms employed in our
formulation. The LGVI, as presented in this paper, preserves the
conserved quantities of the full body dynamics as well as the
orthogonal structure of the rotation matrices. The obtained
simulation results exhibit good stability properties for the
invariants of motion for exponentially long times. The computational
load is further minimized since the LGVI requires one force and
torque evaluation per integration step for second order accuracy.

Subsequent sections of this paper are organized as follows.
Algorithms for the mutual gravitational forces and moments
computation for polyhedral body models are presented in section
\ref{sect:polymutual}. Some description of the full two body problem
and its mathematical properties and the continuous relative
equations of motion are given in section \ref{sect:continuousEOM}.
The Lie Group Variational Integrator for computing the dynamics of
the full body problem is given in section \ref{sect:LGVI}. Selected
simulation results are given in section \ref{sect:results} for two
octahedral rigid bodies in several scenarios, for comparison with
other integration methods and to show the qualities of the results
thus obtained. Conclusions drawn from the results for these
scenarios about the accuracy and suitability of our methods, and
about computational burden, are given in the last section.

\section{Polyhedral Mutual Gravitational Potential, Forces, and Moments}\label{sect:polymutual}

In this section, we present computational algorithms for determining
the mutual gravitational potential, forces, and moments given
polyhedral models of each of the two rigid bodies. The algorithm to
compute the mutual potential relatively efficiently with such a
modelling approach is outlined in Ref. \citen{Werner-CMDA-2005}. The
methods outlined in Ref. \citen{Fahnestock-AASastrospec-2005}, and
in Ref. \citen{Fahnestock-CMDA-2006} with some corrections and improvement of notation, can then be used to compute the gradients of the mutual potential (i.e. forces and moments). A more detailed description and derivation of materials in this section, including extensions to computing the forces and moments for inertial equations of motion, can be found in Ref. \citen{Fahnestock-CMDA-2006}.

\subsection{Context for Methodology}\label{subsect:polymutualdiscussion}

In Ref. \citen{Werner-CMDA-2005}, a uniform density or homogenous
mass distribution within the entire volume of each rigid body is
assumed. This may be a realistic assumption for binary asteroids
based on some empirical data\cite{Yeomans-Science-2000}. However, the
underlying method \emph{does} also allow for approximating arbitrary
density variations, as required for greater asteroid modelling
fidelity and for systems in which one or more bodies is a
spacecraft. This is mentioned without proper explanation in Ref.
\citen{Fahnestock-CMDA-2006}, so we will be more specific here.
Rather than representing each body with just a single polyhedron
whose triangular faces form the body's outer surface, as in prior
work, we can represent each body by multiple polyhedra partially or
fully nested inside of one another. Each polyhedron has a closed
surface consisting of triangular faces, and each triangular face is
defined by three vertices. A tetrahedron is formed by these three
vertices plus the body centroid. This tetrahedron is referred to
hereafter as a simplex. We require each simplex to have a constant
density, but different densities can be assigned to each simplex.
This allows for density variation over the two angular spherical
coordinates within each polyhedron, with resolution determined by
the size of the faces. The possible nesting of multiple polyhedra
allows for density variation over the radial spherical coordinate
within the body volume as well. It should be noted though that to
make practical use of the method herein, the coordinates with
respect to the body centroid of the three non-centroid vertices of
each simplex must be known. While this vertex coordinate information
is known throughout for spacecraft bodies, and can be derived for
exterior faces of natural bodies from shape data obtained via
optical and radar observations or LIDAR surveying, it is difficult
to specify for interior faces within natural bodies without specific
knowledge of internal mass distribution.

The underlying principle of what follows is that the evaluation of the mutual potential's double volume integral over both bodies is equivalent to a global sum of the results of evaluating that double volume integral over each possible pairing of simplices, one being drawn from each body. For each such pairing, the contribution to the mutual potential is given by a Legendre-series expansion. Successive terms of this expansion are linear combinations of other terms which factor into symmetric tensors of increasing rank that are independent of relative position and attitude, and other tensors that depend on the relative position and attitude. The gradients of the mutual potential according to this formulation make use of efficient tensor differentiation rules and the chain rule, resulting in a series expansion for the forces and moments that converges for all points exterior to all polyhedra representing the bodies, and hence for all points exterior to the bodies.

It is useful to note here that the number of terms kept in the
Legendre-series expansion determines the order, in the inverse of
the distance between body centroids, of the errors in the forces and
the moments. When the rigid bodies have small separation distances,
more terms are required to maintain the same error levels. A
possibility for adaptivity during simulation is to adjust the
number of terms used in the series expansion for the force and
moment computations, depending on the body separation distance.
Another possibility is to refine one or more of the polyhedra
representing each body as separation distance decreases. However,
this must be done carefully, in a manner which involves no change to
the relevant conserved quantities and no discontinuity in the motion
states.

It is also worth mentioning that one can consider our approach as a
special case of a more general problem, not treated in this paper,
in which the two bodies are no longer rigid but fragmented, i.e.
they are so-called ``rubble-piles". The rubble-pile model of
asteroids is being supported by an increasing body of evidence from
recent asteroid exploration missions\cite{Fujiwara-Science-2006}. A
number of researchers have represented such rubble-pile asteroids as
collections of nonintersecting spheres held together under
self-gravity\cite{Richardson-Icarus-1998}, not to be confused with
representing a rigid body by filling its volume with spheres fixed
with respect to one another. One could just as well represent
rubble-pile asteroids as collections of ellipsoids or polyhedra of
arbitrary size and shape held together under self-gravity, rather
than smooth spheres, as in Refs. \citen{Roig-Icarus-2003} and
\citen{Korycansky-ASS-2004} respectively. This allows for filling a
body's volume more fully or at different porosities than are
possible with spheres (even zero initial porosity with polyhedra).
With polyhedra, different levels of rigidity within a body can then
be considered simply by initial combination of smaller polyhedra
with coincident or adjacent faces into larger polyhedra. The full
two rigid body problem that is the subject of this paper can then be
viewed from a different perspective as the limiting case of that
process of combination within two separated bodies. This is similar
to using a tree-code method\cite{Duan-JCC-2001}, but only using the
root- or top-level cell (or group) within each body. In other cases
employing multiple non-intersecting polyhedra within each body, the
gravitational force and moment couples between every possible pair
of polyhedra in the binary system can still be accurately evaluated
using the methodology of this section, considered independently from
the rest of this paper.

\subsection{Mutual Gravitational Potential}\label{subsect:polymutualpot}

We label the two polyhedral rigid bodies \B{1} and \B{2}%
\nomenclature[$\B{1}$]{$\B{1}\,,\,\B{2}$}{Labels for the two
polyhedral rigid bodies}
. Consistent with the earlier discussion and definition of the polyhedral modelling approach, let body \B{1} be divided into a set of simplices indexed by $a$%
\nomenclature[$a$]{$a$}{Simplex in \B{1}}%
\nomenclature[b$a$]{$_a$}{Denotes ``for simplex $a$ in \B{1}'' }
 and let body \B{2} be divided into a set of simplices indexed by $b$%
\nomenclature[$b$]{$b$}{Simplex in \B{2}}%
\nomenclature[b$b$]{$_b$}{Denotes ``for simplex $b$ in \B{2}'' } .
Evaluating the double volume integrals over \B{1} and \B{2} is
equivalent to the double summation over all $a$ and over all $b$ of
the result of evaluating the double volume integrals over each
simplex combination ($a$,$b$). This is shown in the following
expression for the mutual potential\cite{Werner-CMDA-2005}. Note
that at this point, we make use of tensor notation and the Einstein
convention of summation over repeated indices:
\begin{multline} \label{eqn:mutpotmain}
U = -G \sum_{a\in \B{1}} \sum_{b\in \B{2}} \rho_a T_a \rho_b T_b
\left\{ \left[\frac{\mathbf Q}{r}\right] + \left[ -\frac{\mathbf Q_i
\mathbf w^i}{r^3} \right] + \left[-\frac{\mathbf Q_{ij} \mathbf
r^{ij}}{2 r^3} + \frac{3 \mathbf Q_{ij} \mathbf w^i \mathbf w^j}{2
r^5}\right] \right. \\ \left.
\;\;\;\;\;\;\;\;\;\;\;\;\;\;\;\;\;\;\;\;\;\;\;\;\;\;\;\;\;\;\;\;\;\;\;
+ \left[ \frac{3 \mathbf Q_{ijk} \mathbf r^{ij} \mathbf w^k}{2 r^5} - \frac{5 \mathbf Q_{ijk} \mathbf w^i \mathbf w^j \mathbf w^k}{2 r^7} \right] + \ldots \right\}%
\nomenclature[$U$]{$U$}{Mutual gravitational potential}%
\nomenclature[$G$]{$G$}{Universal gravitational constant}%
\nomenclature[$\rho$]{$\rho$}{Density, kg/m$^3$}%
\nomenclature[$T$]{$T$}{Jacobian determinant of the matrix containing coordinates of non-centroid vertices of a simplex}%
\nomenclature[$r$]{$r$}{Scalar magnitude of the relative position vector between body centroids, $r=\norm{X}$}%
\nomenclature[a1$Q$]{$\mathbf Q$}{Tensor that is symmetric along every dimension, in which each element is a rational number, with form illustrated in Ref. \citen{Werner-CMDA-2005}}%
\nomenclature[b$i$]{$_i\,_j\,_k\,_p$}{Tensor indices eliminated by summation in left hand side of equations}%
\nomenclature[t$i$]{$^i\,^j\,^k\,^p$}{Tensor indices eliminated by summation in left hand side of equations}%
\nomenclature[b$\phi$]{$_{\phi}\,_{\theta}$}{Tensor indices not eliminated by summation}%
\nomenclature[t$\phi$]{$^{\phi}\,^{\theta}$}{Tensor indices not
eliminated by summation}
\end{multline}
The scalars $T_a$ and $T_b$ and the tensors $\mathbf Q$ of increasing rank are all independent of both relative position between centroids and relative attitude between the bodies, so they can be computed before any dynamic simulation. However the vector $\mathbf w$%
\nomenclature[a2$w$]{$\mathbf w$}{Rank 1 tensor, $\in\Re^6$,
dependent on relative attitude and relative position}
 is dependent on relative position, and both the vector $\mathbf w$ and the matrix $\mathbf r$%
\nomenclature[a2$\mathbf r$]{$\mathbf r$}{Rank 2 tensor,
$\in\Re^{6\times 6}$, dependent on relative attitude but not
relative position}
 are dependent on relative attitude through a matrix $\mathbf v$%
\nomenclature[a2$v$]{$\mathbf v$}{Rank 2 tensor, $\in\Re^{3\times
6}$, containing coordinates of non-centroid vertices of both
simplices in a simplex pairing, expressed in any desired coordinate frame}
, where
\begin{align}
\mathbf w^i \;=\; \mathbf v^i_j X^j ~~~~~~~,~~~~~~~ \mathbf r^{ij} \;=\; \mathbf v^i_p \mathbf v^j_p . \nonumber%
\nomenclature[$X$]{$X$}{Relative position vector between body centroids expressed in the frame fixed to $\B{2}$}%
\nomenclature[t$T$]{$^T$}{Denotes matrix transpose}
\end{align}
The $\mathbf Q$'s are defined in Ref. \citen{Werner-CMDA-2005}, in which they are written out explicitly up to the third rank for illustration of their form. The series within the braces in Eq.~\refeqn{mutpotmain} is infinite but sufficient accuracy seems to be obtained with just the first several terms in square brackets. We denote these bracketed terms as scalars $\hat{U}_0$, $\hat{U}_1$, $\hat{U}_2$,%
\nomenclature[$Uo$]{$\hat{U}_{\#}$}{The $\#$'th Legendre series term
of increasing order within mutual gravitational potential equation}
 and so on.

\subsection{Mutual Gravitational Forces and Moments}\label{subsect:polymutualpotgrads}

The gravitational force existing between the two bodies is given by
the partial derivative of the mutual potential with respect to the
relative position vector $X$. To obtain this we first derive simple
tensor differentiation rules~\cite{Fahnestock-CMDA-2006}:
\begin{align}\label{eqn:forcelemmas}
\frac{\partial r}{\partial X_{\theta}} = \frac{X_{\theta}}{r}
~~~~~~~~,~~~~~~~~ \frac{\partial \mathbf w^i}{\partial X_{\theta}} =
\mathbf v_{\theta}^i.
\end{align}
These are used in finding each successive $\partial
\hat{U}\,/\,\partial X_{\theta}$ term. The first few of these terms
are illustrated below:
\begin{align*}
\deriv{\hat{U}_0}{X_{\theta}}& =-\frac{\mathbf Q X_{\theta}}{r^3},\\
\deriv{\hat{U}_1}{X_{\theta}}& =\frac{3 \mathbf Q_i X_{\theta} \mathbf w^i}{r^5} \;-\; \frac{\mathbf Q_i \mathbf v_{\theta}^i}{r^3},\\
\deriv{\hat{U_2}}{X_{\theta}}& =\frac{3 \mathbf Q_{ij} \mathbf
r^{ij} X_{\theta}}{2r^5}
    \;-\; \frac{15 \mathbf Q_{ij} X_{\theta} \mathbf w^i\mathbf w^j}{ 2 r^7}
    \;+\; \frac{3 \mathbf Q_{ij} \mathbf w^i\mathbf v_{\theta}^j}{r^5},\\
\deriv{\hat{U_3}}{X_{\theta}}& =
  -\; \frac{15 \mathbf Q_{ijk} \mathbf r^{ij} X_{\theta} \mathbf w^k}{2 r^7}
\;+\; \frac{3 \mathbf Q_{ijk} \mathbf r^{ij} \mathbf
v_{\theta}^{k}}{2 r^5} \;+\; \frac{35 \mathbf Q_{ijk} X_{\theta}
\mathbf w^i \mathbf w^j \mathbf w^k}{2 r^9} \;-\; \frac{15 \mathbf
Q_{ijk} \mathbf w^i \mathbf w^j \mathbf v_{\theta}^{k}}{2 r^7}.
\end{align*}
Such terms are used in the overall expression for the force vector
that is in turn used in the relative equations of motion, namely
(laying aside the tensor notation for the moment)
\begin{align} \label{eqn:positiongradient}
\deriv{U}{X} = -G \sum_{a\in \B{1}} \sum_{b\in \B{2}} \rho_a T_a
\rho_b T_b \left( \frac{\partial \hat{U_0}}{\partial X}\;
+\;\frac{\partial \hat{U_1}}{\partial X}\;+
\;\frac{\partial\hat{U_2}}{\partial X}\;+\;\cdots\right).
\end{align}
The torque or moment between the bodies due to their mutual gravitational interaction is given by an expression requiring the partial derivative of the mutual potential with respect to the relative attitude $R$%
\nomenclature[$R$]{$R$}{Relative attitude rotation matrix, $\in\SO$,
mapping from the frame fixed to \B{1} to the frame fixed to \B{2}} .
Evaluating this requires partial differentiation of one attitude
rotation matrix, an element of $\SO$, with respect to another such
matrix. For this we use the basic rule
\begin{align} \label{eqn:RbyRrule}
\deriv{R^{jk}}{R^{\phi\theta}} \;=\; \delta^{\phi}_j \delta^k_{\theta}%
\nomenclature[$\delta$]{$\delta$}{Rank 2 tensor defined by the
Kronecker delta function}
\end{align}
inside of the expression for the tensor $\mathbf D$%
\nomenclature[a$D$]{$\mathbf D$}{Rank 4 tensor,
$\in\Re^{3\times6\times3\times3}$, the partial derivative of
$\mathbf v$ with respect to a rotation matrix} . The details of this
tensor and an alternate rule are discussed in Ref.
\citen{Fahnestock-CMDA-2006}. It is in turn used in the two tensor
differentiation rules~\cite{Fahnestock-CMDA-2006}:
\begin{align} \label{eqn:torquelemmas}
\frac{\partial \mathbf w^i}{\partial R^{\phi \theta}} = X^j \mathbf
D_{j \theta}^{\phi i} ~~~~~~~~,~~~~~~~~ \frac{\partial \mathbf
r^{ij}}{\partial R^{\phi \theta}} = 2\,\mathbf v_p^i \mathbf D_{p
\theta}^{\phi j}
\end{align}
These are used in finding each successive $\partial
\hat{U}\,/\,\partial R^{\phi \theta}$ term. There is no attitude
dependence in $\hat{U}_0$, and then the next few partials are:
\begin{align*}
\frac{\partial \hat{U}_1}{\partial R^{\phi \theta}} & =
 -\frac{\mathbf Q_i X^j\mathbf D_{j \theta}^{\phi
 i}}{r^3},\\
\frac{\partial \hat{U}_2}{\partial R^{\phi \theta}} & =
-\frac{\mathbf Q_{ij} \mathbf v_p^i \mathbf D_{p \theta}^{\phi
j}}{r^3}
\;+\; \frac{3\,\mathbf Q_{ij} \mathbf w^i X^p \mathbf D_{p \theta}^{\phi j}}{r^5},\\
\frac{\partial \hat{U}_3}{\partial R^{\phi \theta}} & =
 \frac{3\;\mathbf Q_{ijk}}{2\,r^5}\left( 2 \mathbf v_p^i \mathbf D_{p \theta}^{\phi j} \mathbf w^k
 \;+\; \mathbf r^{ij} X^p \mathbf D_{p \theta}^{\phi k}\right)
 - \frac{15\,\mathbf Q_{ijk} \mathbf w^i \mathbf w^j X^p \mathbf D_{p \theta}^{\phi k}}{2\,r^7}.
\end{align*}
Such terms are used in the overall expression
\begin{align} \label{eqn:attitudegradient}
\deriv{U}{R^{\phi \theta}} \;=\; -G \sum_{a\in \B{1}} \sum_{b\in
\B{2}} \rho_a T_a \rho_b T_b \left( \frac{\partial
\hat{U}_1}{\partial R^{\phi \theta}} \;+\; \frac{\partial
\hat{U}_2}{\partial R^{\phi \theta}} \;+\; \frac{\partial
\hat{U}_3}{\partial R^{\phi \theta}} \;+\;\dots\right).
\end{align}
Using this within Eq.~\refeqn{M} below, we can evaluate the moment
due to the mutual gravitational potential.

\section{Full Two Body Dynamics and Continuous Equations of Motion}\label{sect:continuousEOM}

In this section, we describe the full two rigid body problem, and we
present the continuous relative equations of motion (EOM) and the
conserved quantities for the full two rigid body dynamics.
Continuous EOM for the full two rigid body problem are given in Ref.
\citen{Maciejewski-CMDA-1995}, and they are formally derived in the
context of Lagrangian mechanics in Ref. \citen{Lee-CMAME-2005}. A
more detailed description and derivation of materials in this
section, including derivation of inertial EOM in addition to
relative EOM, can be found in the latter reference.

The physical configuration of the full two rigid body problem is described in terms of the position vector of each rigid body in an inertial reference frame%
\nomenclature[$x_1$]{$x_1$}{Position vector of $\B{1}$'s centroid expressed in the inertial reference frame}%
\nomenclature[$X_1$]{$X_1$}{Position vector of $\B{1}$'s centroid expressed in the frame fixed to $\B{1}$}%
\nomenclature[$x_2$]{$x_2$}{Position vector of $\B{2}$'s centroid expressed in the inertial reference frame}%
\nomenclature[$X_2$]{$X_2$}{Position vector of $\B{2}$'s centroid expressed in the frame fixed to $\B{2}$}%
, an element of $\Re^3$, and the attitude of each rigid body with respect to that frame.%
\nomenclature[$R_1$]{$R_1$}{Rotation matrix, $\in\SO$, mapping from the frame fixed to $\B{1}$ to the inertial reference frame}%
\nomenclature[$R_2$]{$R_2$}{Rotation matrix, $\in\SO$, mapping from
the frame fixed to $\B{2}$ to the inertial reference frame}
 This attitude is represented by a rotation matrix, which is a $3\times 3$ orthogonal matrix with determinant $+1$, and hence is an element of the Lie group $\SO$. The configuration space of each rigid body is therefore a Lie group referred to as the Euclidean Lie group, $\SE=\Re^3$\textcircled{s}$\SO$ (with \textcircled{s} meaning the semi-direct product). The motion of the two rigid bodies depends only on the relative positions and the relative attitudes of the bodies. This is a consequence of the property that the mutual gravitational potential depends only on these relative variables. Thus, the Lagrangian of the two rigid bodies is invariant under a group action of $\SE$, and it is natural to reduce the EOM for the full body problem by writing them in one of the body fixed frames. Without loss of generality, we present the EOM in the frame fixed to $\B{2}$. Then the configuration space of the reduced system is $\SE$.

We start with the relative variables
\begin{align}
X & = R_2^T(x_1-x_2),\label{eqn:X} \\
R & = R_2^TR_1,\label{eqn:R}
\end{align}
where $X\in\Re^3$ is the relative position of $\B{1}$ with respect
to $\B{2}$ expressed in the frame fixed to $\B{2}$, and $R\in\SO$ is
the relative attitude of $\B{1}$ with respect to $\B{2}$. The
corresponding linear and angular velocities are
\begin{align}
V & = R_2^T(v_1-v_2),\\
\Omega & = R\Omega_1,%
\nomenclature[$v_1$]{$v_1$}{Velocity vector of $\B{1}$'s centroid expressed in the inertial reference frame}%
\nomenclature[$V_1$]{$V_1$}{Velocity vector of $\B{1}$'s centroid expressed in the frame fixed to $\B{1}$}%
\nomenclature[$v_2$]{$v_2$}{Velocity vector of $\B{2}$'s centroid expressed in the inertial reference frame}%
\nomenclature[$V_2$]{$V_2$}{Velocity vector of $\B{2}$'s centroid expressed in the frame fixed to $\B{2}$}%
\nomenclature[$V$]{$V$}{Relative velocity vector of $\B{1}$ with respect to $\B{2}$ expressed in the frame fixed to $\B{2}$}%
\nomenclature[$\Omega_1$]{$\Omega_1$}{Angular velocity vector of $\B{1}$ expressed in the frame fixed to $\B{1}$}%
\nomenclature[$\Omega_2$]{$\Omega_2$}{Angular velocity vector of $\B{2}$ expressed in the frame fixed to $\B{2}$}%
\nomenclature[$\Omega$]{$\Omega$}{Angular velocity vector of $\B{1}$
expressed in the frame fixed to $\B{2}$}
\end{align}
where $V\in\Re^3$ represents the relative velocity of $\B{1}$ with
respect to $\B{2}$ expressed in the frame fixed to $\B{2}$, and
$\Omega\in\Re^3$ is the angular velocity of $\B{1}$ expressed in the
frame fixed to $\B{2}$. The standard and nonstandard moment of
inertia matrices (see Ref. \citen{Lee-CMAME-2005} for the
distinction) of $\B{1}$ can also be expressed in the frame fixed to
$\B{2}$ according to
$J_R=RJ_1R^T$%
\nomenclature[$J_R$]{$J_R$}{Standard moment of inertia matrix of $\B{1}$ expressed in frame fixed to $\B{2}$}%
\nomenclature[$J_1$]{$J_1$}{Standard moment of inertia matrix of $\B{1}$ expressed in its own frame}%
\nomenclature[$J_2$]{$J_2$}{Standard moment of inertia matrix of $\B{2}$ expressed in its own frame}%
\nomenclature[$J_{dR}$]{$J_{dR}$}{Non-standard moment of inertia matrix of $\B{1}$ expressed in frame fixed to $\B{2}$}%
\nomenclature[$J_{d_1}$]{$J_{d_1}$}{Non-standard moment of inertia matrix of $\B{1}$ expressed in its own frame}%
\nomenclature[$J_{d_2}$]{$J_{d_2}$}{Non-standard moment of inertia
matrix of $\B{2}$ expressed in its own frame}
 and $J_{dR}=RJ_{d_1}R^T$. Note that $J_R$ and $J_{dR}$ are not constant matrices.

Assuming that the total linear momentum of the whole system is zero,
the kinetic energy of the system is expressed in terms of the
relative variables as
\begin{align}
KE & =\frac{1}{2}m\norm{V}^2 +\frac{1}{2}\tr{S(\Omega)J_{dR}S(\Omega)^T} +\frac{1}{2}\tr{S(\Omega_2)J_{d_2}S(\Omega_2)^T},%
\nomenclature[$KE$]{$KE$}{Kinetic energy of system}%
\nomenclature[$S(.)$]{$S(.)$}{Denotes the cross-product operation matrix.}
\end{align}
where $m=\frac{m_1m_2}{m_1+m_2}\in\Re$%
\nomenclature[$m$]{$m$}{Scalar mass parameter for the system}%
\nomenclature[$m_1$]{$m_1$}{Mass of body $\B{1}$}%
\nomenclature[$m_1$]{$m_1$}{Mass of body $\B{1}$}%
\nomenclature[$m_2$]{$m_2$}{Mass of body $\B{2}$} , and
$S(\cdot):\Re^3\mapsto \so$ is an isomorphism between $\Re^3$ and
skew-symmetric matrices, defined such that $S(x)y=x\times y$ for any
$x,y\in\Re^3$.
 We also have from Eq.~\refeqn{mutpotmain} the potential energy expressed in terms of the relative variables, that is to say $PE=U(X,R)$%
\nomenclature[$PE$]{$PE$}{Potential energy of system}, recalling that in Eq.~\refeqn{mutpotmain} $r$ and $\mathbf w$ are functions of $X$ and both $\mathbf w$ and $\mathbf r$ are functions of $R$. The following continuous EOM of the full two rigid body problem in relative coordinates can be obtained by taking variations of the reduced Lagrangian, $L=KE-PE$,%
\nomenclature[$L$]{$L$}{Reduced Lagrangian for system} as follows:
\begin{gather}
\dot{X}+\Omega_2 \times X=V,\\
\dot{R}=S(\Omega)R-S(\Omega_2)R,\label{eqn:Rdot}\\
m\dot{V}+m\Omega_2\times V = -\deriv{U}{X},\label{eqn:Vdot}\\
\dot{(J_R\Omega)}+\Omega_2\times J_R\Omega=-M,\\
J_2\dot{\Omega}_2+\Omega_2\times J_2\Omega_2=X \times
\deriv{U}{X}+M,\label{eqn:Pi2dot}
\end{gather}
where the vector $J_R\Omega=RJ_1\Omega_1\in\Re^3$ is the angular
momentum of the first body expressed in the second body fixed frame.
The moment
due to the gravity potential, $M\in\Re^3$%
\nomenclature[$M$]{$M$}{Moment due to mutual gravitational
potential}, is determined by the following relationship
\begin{align*}
S(M)=\deriv{U}{R}R^T-R \deriv{U}{R}^T,
\end{align*}
or more explicitly,
\begin{align}\label{eqn:M}
M = R_{c1} \times \left[\deriv{U}{R}\right]_{c1} + R_{c2} \times \left[\deriv{U}{R}\right]_{c2} + R_{c3} \times \left[\deriv{U}{R}\right]_{c3}.%
\nomenclature[b$c1$]{$_{c1}\,_{c2}\,_{c3}$}{Denotes the first,
second, third column vectors of a matrix, respectively}
\end{align}
These equations are also given in different notation in Ref.
\citen{Fahnestock-CMDA-2006}. The total energy, the total linear
momentum expressed in the inertial frame, and the total angular
momentum expressed in the inertial frame are conserved quantities,
given by
\begin{align}
E & = KE + PE\\
\gamma_T & = m_1 v_1 + m_2 v_2,\\
\pi_T & = x_1 \times m_1 v_1 + R_1 J_1\Omega_1 + x_2 \times m_2 v_2 + R_2 J_2\Omega_2.%
\nomenclature[$E$]{$E$}{Total energy of system}%
\nomenclature[$\gamma_T$]{$\gamma_T$}{Total linear momentum of system expressed in inertial reference frame}%
\nomenclature[$\pi_T$]{$\pi_T$}{Total angular momentum of system
expressed in inertial reference frame}
\end{align}

\section{Discrete Equations of Motion: the Lie Group Variational Integrator}\label{sect:LGVI}

This section presents the LGVI, essentially a statement of discrete relative EOM in Hamiltonian form, as opposed to continuous relative EOM in Lagrangian form in the previous section. Again, a detailed description and derivation of materials in this section, including development of inertial discrete EOM in addition to these relative discrete EOM, can be found in Ref. \citen{Lee-CMAME-2005}. Here we simply state the results.
\nomenclature[b$n$]{$_n$}{Second level subscript, denotes the value of variables at simulation time $t=n h+t_0$}%
\nomenclature[$t$]{$t$}{Simulation time}%
\nomenclature[$t_0$]{$t_0$}{Time at start of simulation}%
\nomenclature[$h$]{$h$}{Integration step size}
\nomenclature[$N$]{$N$}{Number of time-steps of length $h$ to go
from the initial time to the final time in simulation}
\nomenclature[$t_f$]{$t_f$}{Time at end of simulation}

\subsection{Discrete Equations of Motion}\label{subsect:discreteEOM}

Since the dynamics of the full two rigid bodies arise from
Lagrangian or Hamiltonian mechanics, they are characterized by
symplectic, momentum and energy preserving properties. These
geometric features determine the qualitative behavior of the
dynamics, and they can serve as a basis for further theoretical
study of the full two rigid body problem. Furthermore, the
configuration space has a group structure denoted by $\SE$.

However, general numerical integration methods, including the
popular Runge-Kutta methods, fail to preserve these geometric
characteristics. General integration methods are obtained by
approximating continuous EOM by directly discretizing them with
respect to time. With each integration step, the updates involve
additive operations, so that the  underlying Lie group structure is
not necessarily preserved as time progresses. This is caused by the
fact that the Euclidean Lie group is not closed under addition.

For example, if we use a Runge-Kutta method for numerical
integration of \refeqn{Rdot}, then the rotation matrices drift from
the orthogonal rotation group, $\SO$; the quantity $R^T R$ drifts
from the identity matrix. Then the attitudes of the rigid bodies
cannot be determined accurately, resulting in significant errors in
the gravitational force and moment computations that depend on the
attitude, and consequently errors in the entire simulation. It is
often proposed to parameterize \refeqn{Rdot} by Euler angles or unit
quaternions. However, Euler angles are not global expressions of the
attitude since they have associated singularities. Unit quaternions
do not exhibit singularities, but are constrained to lie on the unit
three-sphere $\mathbb{S}^3$, and general numerical integration
methods do not preserve the unit length constraint. Therefore,
quaternions lead to the same numerical drift problem. Re-normalizing
the quaternion vector at each step tends to break the conservation
properties. Furthermore, unit quaternions double cover $\SO$, so
that there are inevitable ambiguities in expressing the attitude.

In contrast, the Lie Group Variational Integrator has desirable
properties such as symplecticity, momentum preservation, and good
energy stability for exponentially long time periods. It also
preserves the Euclidian Lie group structure without the use of local
charts, reprojection, or constraints. The LGVI is obtained by
discretizing Hamilton's principle; the velocity phase space of the
continuous Lagrangian is replaced by discrete variables, and a
discrete Lagrangian is chosen such that it approximates a segment of
the action integral. Taking the variation of the resulting action
sum, we obtain discrete EOM referred to as a variational integrator.
Since the discrete variables are updated by Lie group operations, the
group structure is preserved. Here we present the resulting discrete
EOM as follows; the detailed development can be found in Ref.
\citen{Lee-CMAME-2005}.

\begin{gather}
X_{_{n+1}}=F_{2_{n}}^T\parenth{X_{_{n}}+hV_{_n}-\frac{h^2}{2m}\deriv{U_{_n}}{X_{_{n}}}},
\label{eqn:updateX}\\
h S\parenth{J_{R_n}\Omega_{_n}-\frac{h}{2}M_{_{n}}}=F_{_{n}}J_{dR_n}-J_{dR_n}F_{_{n}}^T,\label{eqn:findF}\\
h
S\parenth{J_2\Omega_{2_{n}}+\frac{h}{2}X_{_{n}}\times\deriv{U_{_n}}{X_{_{n}}}+\frac{h}{2}M_{_{n}}}
=F_{2_{n}}J_{d_2}-J_{d_2}F_{2_{n}}^T,\label{eqn:findF2}\\
R_{_{n+1}}=F_{2_{n}}^TF_{_{n}}R_{_{n}},\label{eqn:updateR}\\
V_{_{n+1}}=F_{2_{n}}^T\parenth{V_{_{n}}-\frac{h}{2m}\deriv{U_{_n}}{X_{_{n}}}}
-\frac{h}{2m}\deriv{U_{_{n+1}}}{X_{_{n+1}}},
\label{eqn:updateV}\\
\nomenclature[$F$]{$F$}{Update matrix, $\in\SO$, for relative attitude rotation matrix $R$}%
\nomenclature[$F_2$]{$F_2$}{Update matrix, $\in\SO$, for rotation
matrix $R_2$} J_{R_{n+1}}\Omega_{_{n+1}}=F_{2_n}^T
\parenth{J_{R_n}\Omega_{_n}-\frac{h}{2}M_{_n}}-\frac{h}{2}M_{_{n+1}},
\label{eqn:updatePi}\\
J_2\Omega_{2_{n+1}}=F_{2_n}^T\parenth{J_2\Omega_{2_n}+\frac{h}{2}X_{_n}\times\deriv{U_{_n}}{X_{_n}}
+\frac{h}{2}M_{_n}}+\frac{h}{2}X_{_{n+1}}\times\deriv{U_{_{n+1}}}{X_{_{n+1}}}+\frac{h}{2}M_{_{n+1}}.
\label{eqn:updatePi2}
\end{gather}
To propagate these equations, we start with a set of initial states,
$(X_{_0},V_{_0},R_{_0},\Omega_{_0},\Omega_{2_0})$, and perform one
initial evaluation of the mutual potential gradients, obtaining
$\partial U_{_0} / \partial X_{_0}$ and $M_{_0}$ with Eqs.
\refeqn{positiongradient} and \refeqn{M}. We then find $X_{_1}$ from
Eq.~\refeqn{updateX}. Solving the implicit equations \refeqn{findF}
and \refeqn{findF2} yields the matrix-multiplication update matrices
$F_{_0}$ and $F_{2_0}$ for the attitude rotation matrices, and
$R_{_1}$ follows from Eq.~\refeqn{updateR}. After that, we use
$X_{_1}$ and $R_{_1}$ in a new evaluation of the mutual potential
gradients. We then compute $V_{_1}$, $\Omega_{_1}$, and
$\Omega_{2_1}$ from equations \refeqn{updateV}, \refeqn{updatePi}
and \refeqn{updatePi2}, respectively. This yields a discrete map
$(X_{_0},V_{_0},R_{_0},\Omega_{_0},\Omega_{2_0})\mapsto
(X_{_1},V_{_1},R_{_1},\Omega_{_1},\Omega_{2_1})$, and this process
can be repeated for each time step. Note that only one new evaluation of the potential
gradients is required per time step. The discrete trajectory in
reduced variables can be used to reconstruct the inertial motion of
the bodies. Either concurrently with that propagation or later after
completion of it, through storing values, we can use the gradient
$\partial U /
\partial X$, the relative attitude $R$, and the update matrix
$F_{2}$ with these equations:
\begin{gather}
x_{2_{n+1}}=x_{2_{n}}+h v_{2_{n}} + \frac{h^2}{2 m_2} R_{_n}
\deriv{U_{_n}}{X_{_n}},
\label{eqn:updatex2}\\
v_{2_{n+1}}=v_{2_{n}} + \frac{h}{2 m_2} R_{_n}
\deriv{U_{_n}}{X_{_n}} + \frac{h}{2 m_2} R_{_{n+1}}
\deriv{U_{_{n+1}}}{X_{_{n+1}}},
\label{eqn:updatev2}\\
R_{2_{n+1}}=R_{2_n} F_{2_n}. \label{eqn:updateR2}
\end{gather}

In the discrete map defined by the LGVI above, the only implicit
parts are Eqs. \refeqn{findF} and \refeqn{findF2}. These two
equations have the same structure, which suggests a specific
computational approach. Using the Rodrigues formula, we rewrite
those equations as an equivalent vector equations, and we solve them
numerically using Newton's iteration. Numerical simulations show
that two or three iterations are sufficient to achieve a tolerance
of $\epsilon=10^{-15}$.

\subsection{Properties of the Lie Group Variational Integrator}\label{subsect:LGVIproperties}

Since the LGVI is obtained by discretizing Hamilton's principle, it
is symplectic and preserves the structure of the configuration
space, \SE, as well as the relevant geometric features of the full
two rigid body problem dynamics represented by the conserved first
integrals of total angular momentum and total energy. The total
energy oscillates around its initial value with small bounds on a
comparatively short timescale, but there is no tendency for the mean
of the oscillation in the total energy to drift (increase or
decrease) from the initial value for exponentially long time. In
contrast, the total energy behavior with general numerical methods
such as the Runge-Kutta schemes tends to drift dramatically over
exponentially long time.

The LGVI preserves the group structure. By using the given
computational approach, the matrices $F_{_n}$ and $F_{2_n}$,
representing the change in relative attitude and attitude of $\B{2}$
over a time step, are guaranteed to be rotation matrices. The group
operation of the Lie group $\SO$ is matrix multiplication. Hence
rotation matrices $R_{_n}$ and $R_{2_n}$ are updated by the group
operation in Eqs. \refeqn{updateR} and \refeqn{updateR2}, so that
they evolve on $\SO$ automatically without constraints or
reprojection. Therefore, the orthogonal structure of the rotation
matrices is preserved, and the attitude of each rigid body is
determined accurately and globally without the need to use local
charts (parameterizations) such as Euler angles or quaternions.

This geometrically exact numerical integration method yields a
highly efficient and accurate computational algorithm, especially
for the full two rigid body problem examined here. In the full two
rigid body problem there is a large burden in computing the mutual
gravitational force and moment for arbitrary bodies, so the number
of force and moment evaluations should be minimized. We have seen
that the LGVI requires only one such evaluation per integration
step, the minimum number of evaluations consistent with the
presented LGVI having second order accuracy (because it is a
self-adjoint method). Within the LGVI, two implicit equations must
be solved at each time step to determine the matrix-multiplication
updates for $R$ and $R_2$. However the LGVI is only weakly implicit
in the sense that the iteration for each implicit equation is
independent of the much more costly gravitational force and moment
computation. The computational load to solve each implicit equation
is comparatively negligible; only two or three iterations are required.
Altogether, the entire method could be considered ``almost
explicit".

The LGVI is a fixed step size integrator, but all of the properties
above are independent of the step size. Consequently, we can achieve
the same level of accuracy while choosing a larger step size as
compared to other numerical integrators of the same order. All of
these features are revealed in the simulation results below.

\section{Numerical Examples}\label{sect:results}

\subsection{Implementation Details}\label{sect:implementation}
Here we describe our specific implementation of the combined
computational methods outlined in this paper. Our codes consist of
pre-processing scripts and post-processing scripts written in the
MATLAB scripting language, and actual executables written in the C
language.

Body models for a natural body or for a spacecraft are originally
provided to us already in the form of ordered vertex and face lists.
Each row of the vertex list contains three numbers for the
coordinates of that vertex's position with respect to a specified
reference frame, preferably with origin at the center of mass of the
body and axes aligned with the body's principle axis. Each row of
the face list contains the three row numbers for the rows of the
vertex list corresponding to the three vertices that form the
corners of that face. The row numbers reference the vertex
information in order, moving counterclockwise around the vertices as
viewed from outside of the body (faces are oriented with
outward-pointing normal vectors according to the right hand rule).
If enough information about density distribution for the body is
known, so that a density number is also assigned to the simplex
associated with each face, those numbers are present in an
additional column in the ordered face list. Otherwise, all simplices
may be given the same density value. In this paper we do not address
the process of generating such formatted body model files from
actual observation data for asteroids or from actual CAD software
models for spacecraft. The face and vertex lists can be manually
generated from scratch for simple arbitrary polyhedron shapes, as
for the octahedra bodies used in the next section. In addition to
the body models themselves, the initial conditions are needed, including the
initial attitude of each body with respect to some common reference
frame, the initial spin axis orientation and spin rate of each body
in that frame, and the initial mutual orbital elements or equivalent
relative translational motion parameters with respect to that frame.

The pre-processing scripts handle a number of preliminary
computations. For each body, the centroid (center of mass) and the
principal axis directions are found, and if not at the origin and
along the elementary unit vectors, respectively, of the frame in
which the vertex coordinates are given, the coordinates of every
vertex are shifted and rotated so that this becomes the case. The
moment of inertia contribution of each simplex in the new frame is
found, and these are summed to obtain the total standard moment of
inertia matrix. Other information such as the surface area, volume,
and mean equivalent radius of each body is also found and reported.
All vertex coordinates and other quantities can then optionally be
nondimensionalized using user-entered length, mass, and time
factors, for better numerical conditioning. Finally, five files are
produced. The first two of these files, one for each body, contains
a re-ordering of the elements of the position vectors, with respect
to the centroid, for the vertices of the simplices. The other three
files, respectively, contain the initial state vector
$(X_{_0},mV_{_0},J_2\Omega_{2_0},J_{R_0}\Omega_{_0},R_{_0},R_{2_0})$,
other system physical data (densities, body volumes, $J_1$, $J_2$,
$m_1$, $m_2$, $m$, nondimensionalized $G$, and nondimensionalization
factors), and the integration parameters (starting and stopping
times and truncation error tolerance, if needed).

A first executable makes use of the Runge-Kutta-Fehlberg 7(8)
integration method, hereafter referred to as RKF7(8) or just RKF, to
propagate the system, after a starting calculation of the
time-invariant $\mathbf{Q}$ tensors. It is noted that this one-time
only calculation of the $\mathbf Q's$ is analogous to finding the
successively higher-order mass moments of a ``normalized" simplex
with vertices at $(0,0,0)$, $(1,0,0)$, $(0,1,0)$, and $(0,0,1)$. It
is also noted that for this high-order scheme, the EOM are evaluated
thirteen times within each integration step, and an evaluation of
the mutual gravitational potential force and moment is required each
of those times. After this the state update is performed and the new
state vector is written to an output file only if the truncation
error is within the tolerance specified. Step size adjustment is
performed with every step. With each state update the mutual
potential itself, force, and moment, while not required for the
propagation of the dynamics, are also evaluated again using the new
state vector and written to another output file. This allows for
checking the total energy conservation and the linear and angular
accelerations. Hence fourteen force and moment evaluations are
needed per integration step.

Another executable uses the LGVI. At the start of this, the
initial state vector from the input file,
$(X_{_0},mV_{_0},J_2\Omega_{2_0},J_{R_0}\Omega_{_0},R_{_0},R_{2_0})$,
is converted to the vector for the discrete mapping, $(X_{_0},V_{_0},R_{_0},\Omega_{_0},\Omega_{2_0})$.
Again, rather than fourteen force and moment evaluations, this LGVI
involves only one such evaluation per integration step, plus the
quick solution of the implicit equations. Step size is fixed, but
this does not present a significant problem for most scenarios
observed in binary asteroids having mutual orbits with relatively
low-eccentricity.

With the aim of completing simulations much faster than in a
single-processor environment, a parallelized version of each
executable above was also written in C with the addition of
MPI. In using the methods presented in this paper, most of the
computation time is associated with the evaluation of the potential
gradients, and that involves performing the same operations for all
of the different simplex combinations, followed by a global sum.
This is well-suited for parallelization. The parallelized version of
each executable is flexible in that any number of nodes or
processors can be specified by the user. Then the process 0 assigns
to each of the other processors the task of calculating the portion
of the potential gradients double summations of Eqs.
\refeqn{positiongradient} and \refeqn{attitudegradient} that arises
from pairing a single simplex in $\B{2}$ with all simplices in
$\B{1}$ in succession. If the number of other processors specified
by the user or found available on the cluster is less than the
number of simplices in $\B{2}$, this is done in rounds until the
portion of the problem matching with every $a$ is obtained. The
parallelized version of each executable has been used on
Myrinet clusters at the Center for Advanced Computing (CAC) at the
University of Michigan and at the Supercomputing and Visualization
Center at NASA's Jet Propulsion Laboratory. Though compiler and user
environment differences produced markedly different capabilities in
each cluster environment, eventually a further two orders of
magnitude reduction in computation time over otherwise identical
single-processor runs was achieved with both the VI and RKF7(8)
schemes. It should be noted however, that for the simulations with
results presented in the following section, the parallel computing
capability was not needed for such small (in number of faces) body
models, and was not utilized.

The MATLAB script post-processing of the output files generates all
desired plots of various dynamic quantities, with optional
capability for animation generation. All pre- and post-processing
steps and scripts are identical regardless of whether the RKF7(8) or
LGVI executable is used, and regardless of whether the
single-processor or parallel version of each is used.

\subsection{Simulation Results}
Simulation results for two octahedral rigid bodies with eight faces
and eight simplices each are given for a variety of scenarios.
Octahedra are used rather than more complex shapes because they are
the simplest polyhedral shapes that manifest the coupled dynamics
behavior desired in all of the scenarios. For greater simplicity,
the octahedra are made symmetric about all axes, although they are
of different sizes. The extents data defining the locations of the corners
of each octahedron are given in Table~\reftab{bodyprops}, as are
various physical parameters of each octahedron including mass and
moment of inertia properties. We present simulation results for four
scenarios, and the initial condition for each scenario is given in
Table~\reftab{ICs}.

\begin{table}
 \begin{center}
  \caption{Properties of octahedral body models used in simulations}\label{tab:bodyprops}
  \begin{tabular}{l|r|r|r|r}
  \multicolumn{1}{c|}{Property} & \multicolumn{2}{c|}{$\B2$} & \multicolumn{2}{c}{$\B1$} \\
\hline
Surface area (m$^2$) & \multicolumn{2}{r|}{8.839} & \multicolumn{2}{r}{2.002} \\
Volume (m$^3$) & \multicolumn{2}{r|}{1.800} & \multicolumn{2}{r}{0.1561} \\
Equiv. radius (m) & \multicolumn{2}{r|}{0.7546} & \multicolumn{2}{r}{0.3340} \\
Mass (kg) & \multicolumn{2}{r|}{4500} & \multicolumn{2}{r}{390.3} \\
Density (kg/m$^3$) & \multicolumn{2}{r|}{2500} & \multicolumn{2}{r}{2500} \\
$I_{xx}$ (kg-m$^2$) & \multicolumn{2}{r|}{1377.0} & \multicolumn{2}{r}{9.24} \\
$I_{yy}$ (kg-m$^2$) & \multicolumn{2}{r|}{814.5} & \multicolumn{2}{r}{42.99} \\
$I_{zz}$ (kg-m$^2$) & \multicolumn{2}{r|}{1462.5} & \multicolumn{2}{r}{44.32} \\
\hline
Extents (m) & min & max & min & max \\
\hline
body frame X & -1.0 & 1.0 & -1.0 & 1.0 \\
body frame Y & -1.5 & 1.5 & -1/$\exp(1)$ & 1/$\exp(1)$ \\
body frame Z & -0.9 & 0.9 & -1/$\pi$ & 1/$\pi$
  \end{tabular}
 \end{center}
\end{table}

\begin{table}[h!]
 \begin{center}
  \begin{threeparttable}
   \caption{Initial Conditions}\label{tab:ICs}
    \begin{tabular}{c|c|c|c}\hline\hline
   Scenario & Attitude$^*$ (deg) & Body spin$^{**}$ (rad/s) & Orbital elements (m,deg) OR state vector (m,m/s)\\\hline
    \multirow{2}*{1} &
    \multirow{2}{1.0in}{(100, 9.8, 175) \\ (160, -5, 165)} &
    \multirow{2}{1.3in}{(0, 0, $5.0\times 10^{-5}$) \\ (0, 0, $9.2\times 10^{-5}$)} &
    \multirow{2}*{($a,e,i,\Omega,\omega,\nu$)=(4m, 0.3, 5$^\circ$, 15$^\circ$, 60$^\circ$, 10$^\circ$)}\\ & & &\\\hline
    \multirow{2}*{2} &
    \multirow{2}{1.0in}{(180, 0, 30)\\ (270, 0, 30)} &
    \multirow{2}{1.3in}{(0, 0, 0.566)\\ (0, 0, -0.566)} &
    \multirow{2}*{$X_0=[0, 6, 0]^T$, $V_0=[-0.006, 0, 0]^T$}\\ & & &\\\hline
    \multirow{2}*{3} &
    \multirow{2}{1.0in}{(-22.6, 5, 180) \\ (50.3, 5, -180)} &
    \multirow{2}{1.3in}{(0, 0, $1.63\times 10^{-4}$) \\ (0, 0, $1.55\times 10^{-4}$)} &
    \multirow{2}*{($a,e,i,\Omega,\omega,\nu$)=(52.9m, 0.942, 5$^\circ$, 0$^\circ$, 88.2$^\circ$, -107.1$^\circ$)}\\ & & &\\\hline
    \multirow{2}*{4} &
    \multirow{2}{1.0in}{(-75, 30, 180)\\ (-75, 30, 180)} &
    \multirow{2}{1.3in}{(0.007, 0.007, 0.05)\\ (-0.003, 0.002, 0.004)} &
    \multirow{2}*{$X_0=[-0.5, 1.8, 1.1]^T$, $V_0=[-0.3, -0.1, 0]^T$}\\ & & &\\\hline\hline
   \end{tabular}
   \begin{tablenotes}
    \item[*] 3-1-3 Euler sequence for $\B1$ (first line) and $\B2$ (second line).
    \item[**] Components of angular velocity of each body expressed in its own body-fixed frame for $\B1$ (first line) and $\B2$ (second line).
   \end{tablenotes}
  \end{threeparttable}
 \end{center}
\end{table}

\paragraph{Scenario 1} The first scenario presented here is that of short-duration simulation of the two octahedra starting from initial conditions matching with a medium eccentricity elliptical mutual orbit.
Both the RKF7(8) and LGVI integrators are used, with the intent of
making a direct comparison between the trajectories of the
configuration variables that result from using each integrator over
a short simulation duration.

\begin{figure}[htb]\centering
 \begin{subfigmatrix}{2}
  \subfigure[Inertial position for $\B2$]{\includegraphics[width=0.45\textwidth]{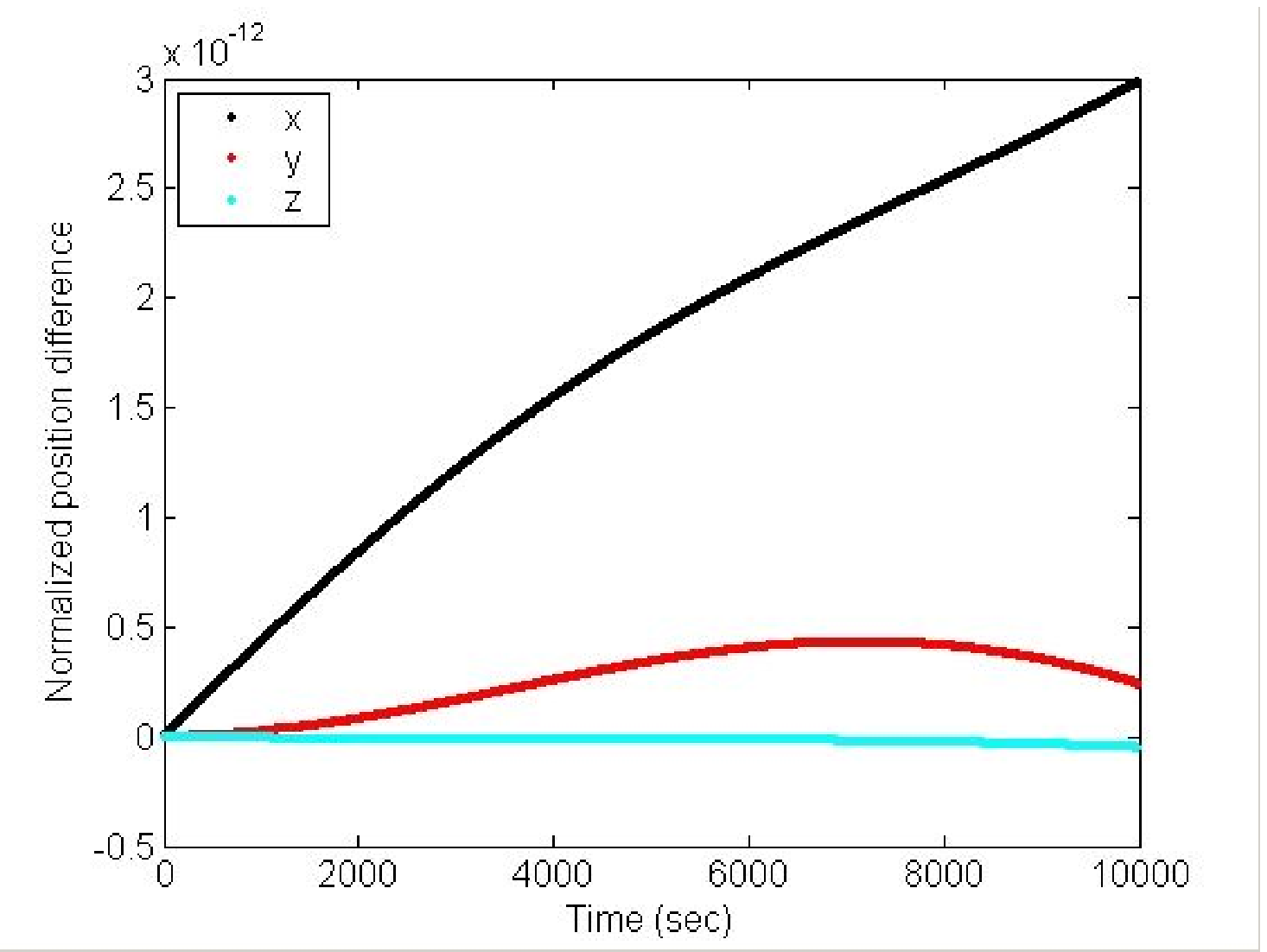} \label{fig:scenario1diff_a}}
  \subfigure[Inertial velocity for $\B2$]{\includegraphics[width=0.45\textwidth]{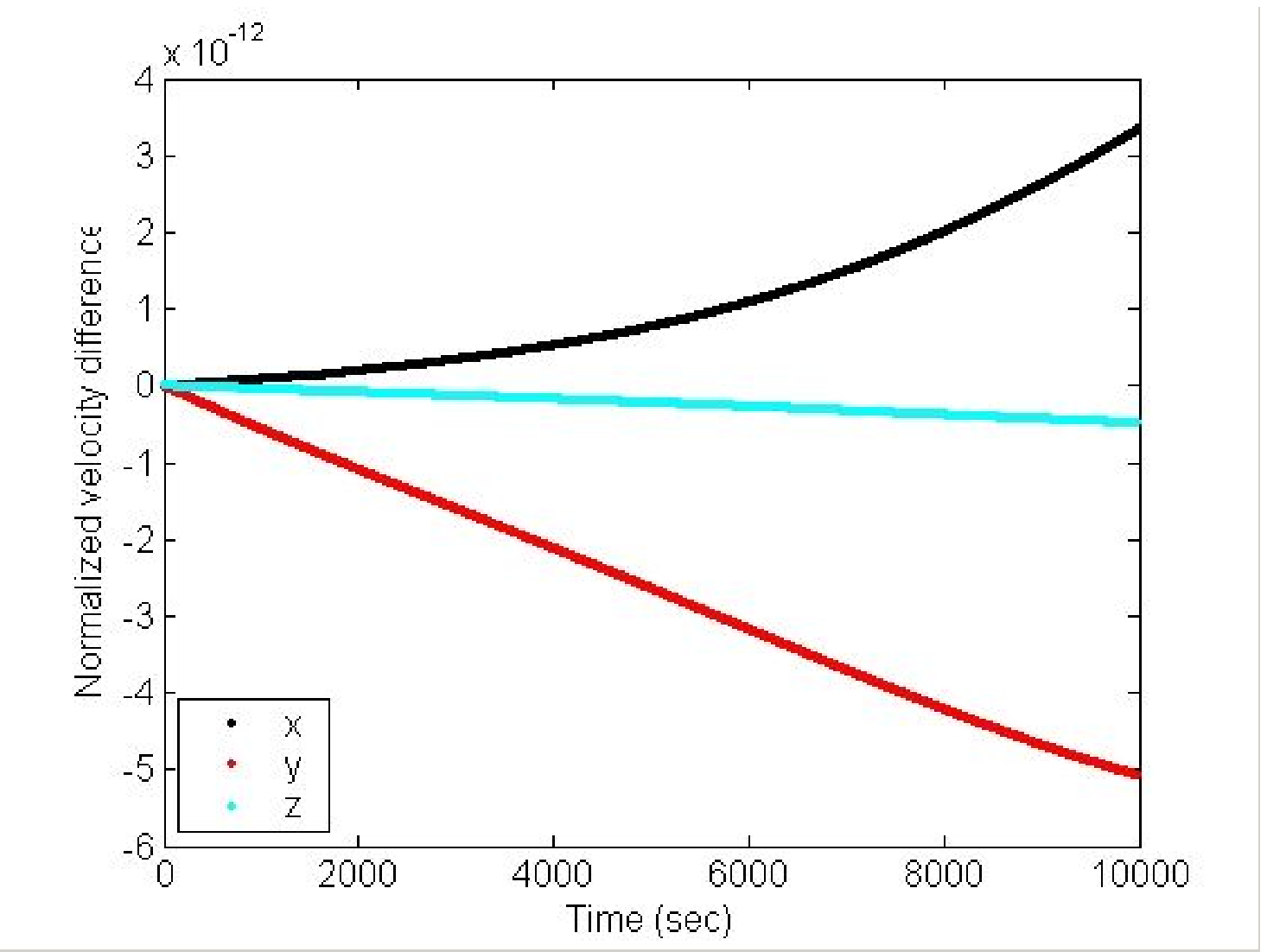} \label{fig:scenario1diff_b}}
  \subfigure[Angular velocity for $\B2$ ]{\includegraphics[width=0.45\textwidth]{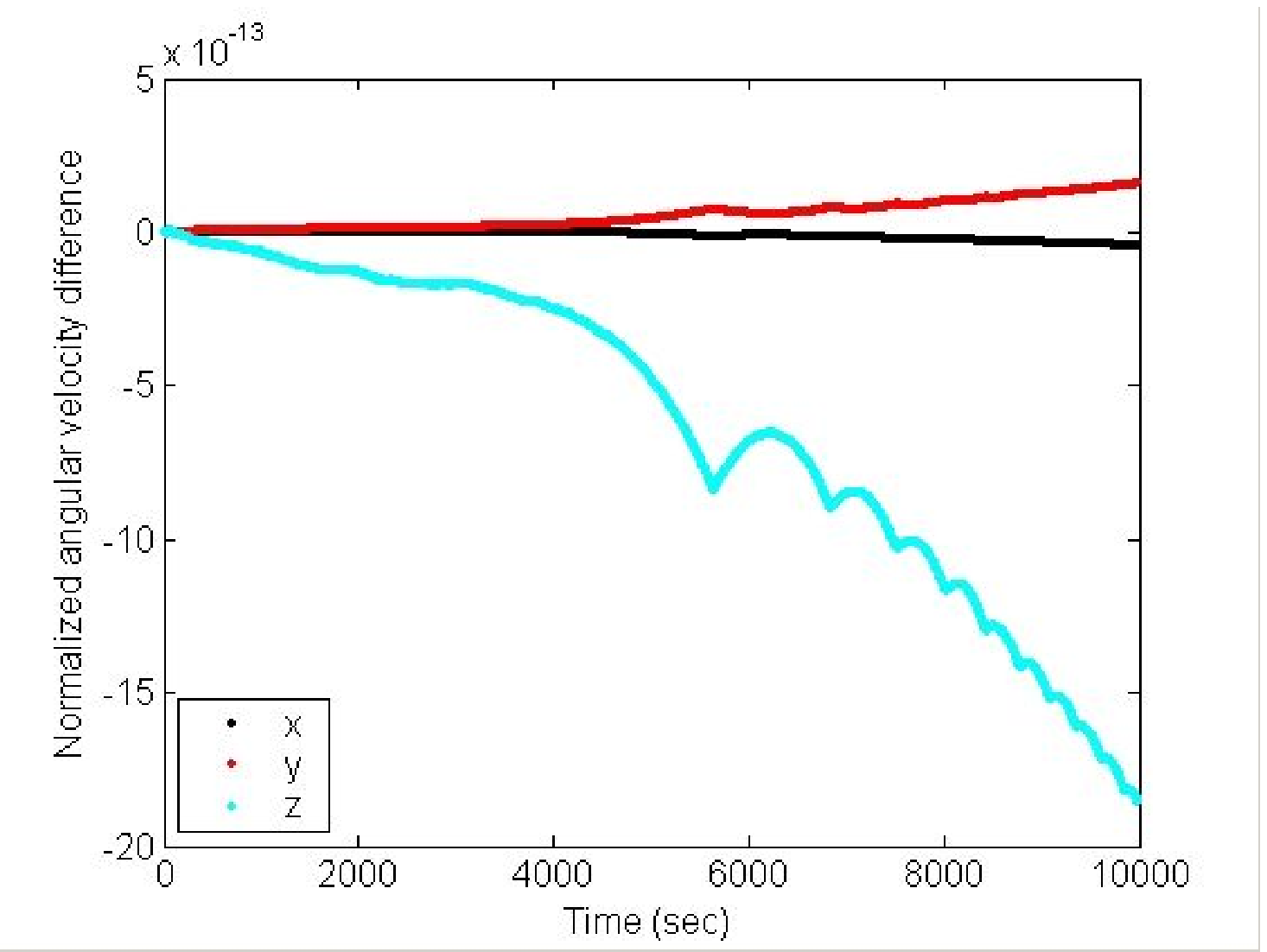} \label{fig:scenario1diff_c}}
  \subfigure[3-1-3 Euler angles for $\B2$]{\includegraphics[width=0.45\textwidth]{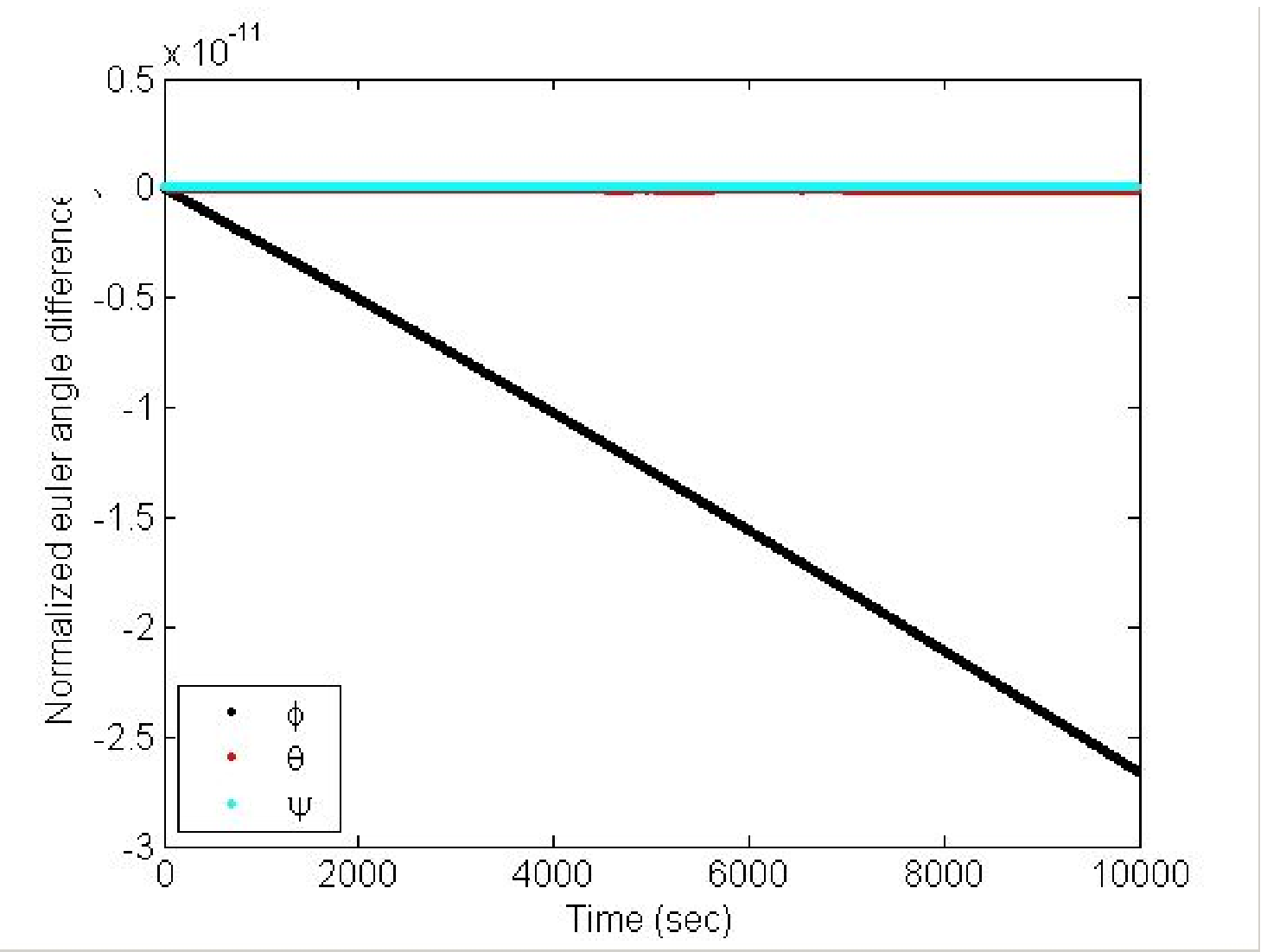} \label{fig:scenario1diff_d}}
\end{subfigmatrix}
\caption{(Scenario 1) Difference between RKF7(8) and LGVI output.
\label{fig:scenario1diff}}
\end{figure}

Figure~\reffig{scenario1diff} shows the difference between the
output of the RKF7(8) and that of the LGVI in components of
reconstructed inertial position, inertial velocity, and body-frame
angular velocity vectors for $\B2$, plus the difference in body
attitude parameters for $\B2$. The corresponding output difference
plots for body $\B1$ look very similar. The differences in vector
components of Figure~\reffig{scenario1diff_a} are normalized by the
system's semi-major axis ($a = 4.0$ m). The differences in vector
components of Figure~\reffig{scenario1diff_b} are normalized by the
equivalent circular velocity ($\sqrt{\mu/a} = 2.856\times10^{-4}$
m/s), and those of Figure~\reffig{scenario1diff_c} are normalized by
the equivalent meanmotion ($\sqrt{\mu/a^3} = 7.141\times10^{-5}$
radians/s). To obtain the results compared here, the  total number
of mutual potential derivatives evaluations and actual running time
using the RKF7(8) routine were 70014 and 494 seconds, respectively,
while the number of such evaluations and actual running time using
the LGVI were 70001 and 539 seconds. Therefore the computational
effort and resources used were roughly the same in each case. All of
these results show that the LGVI can be trusted to produce almost
exactly the same trajectory as a standard RKF7(8) integration
routine over short time scales. As the simulation duration increases
the trajectories from the two integrators begin to diverge. The
behavior of integrals of motion and appropriate error metrics must
then be used to discern which trajectory is to be taken as the
``truth".

\paragraph{Scenario 2} Another scenario is that of propagation from an initial condition with the bodies aligned but possessing relatively large magnitude centroid velocity vectors that are antiparallel and perpendicular to the initial line between centroids. This scenario is simulated both with the LGVI at different step sizes and with the variable step size RKF7(8) at different error tolerances. This allows for a comparison between the integrators of their performance, in terms of the total energy and total angular momentum integrals and the attitude error metric growth vs. computational burden. The results in Table~\reftab{perfcompare} illustrate the general
superiority of the LGVI approach over Runge-Kutta-type approaches.
\begin{table}[htb]
 \begin{center}
  \begin{threeparttable}
  \caption{(Scenario 2) Performance comparison between RKF7(8) and LGVI \label{tab:perfcompare}}   \begin{tabular}{l|cccc|ccc}\hline\hline
Method & $h^*$ & $N_u^{\star}$ & $t_{W}^{\diamond}$ & $\epsilon^{\triangleleft}$
& $\mathrm{E}[|\Delta\mbox{TE}|]^{\dagger \ddag}$ & $\mathrm{E}[\norm{\Delta\pi_T}]^{\dagger \ddag}$ & $\mathrm{E}[\norm{I-R^TR}]^{\dagger}$ \\
\hline
RKF7(8) & $0.236$ & $2368912$ & $23439$ & $10^{-12}$ & $3.901\times10^{-12}$ & $1.493\times10^{-9}$ & $1.151\times 10^{-7}$ \\
RKF7(8) & $0.421$ & $1331414$ & $9102$ & $10^{-10}$ & $1.274\times10^{-10}$ & $2.630\times10^{-7}$ & $1.985\times 10^{-5}$ \\
RKF7(8) & $0.749$ & $747376$ & $5252$ & $10^{-8}$ & $2.284\times10^{-8}$ & $4.620\times10^{-5}$ & $3.173\times 10^{-3}$ \\
\hline
LGVI& $0.0169 $ & $2370000$ & $13511$ & - & $1.698\times 10^{-11}$ & $5.167\times 10^{-10}$ & $2.525\times 10^{-11}$ \\
LGVI& $0.04 $ & $1000000$ & $9920$ & - & $1.928\times 10^{-11}$ & $1.189\times 10^{-10}$ & $2.120\times 10^{-11}$ \\
LGVI& $0.08 $ & $500000$ & $5127$ & - & $9.879\times 10^{-11}$ & $4.139\times 10^{-11}$ & $2.004\times 10^{-12}$ \\
LGVI& $0.4$ & $100000$ & $983$ & - & $2.234\times 10^{-9}$ & $6.266\times 10^{-12}$ & $3.386\times 10^{-14}$ \\
LGVI& $0.8$ & $50000$ & $431$ & - & $9.326\times 10^{-9}$ & $1.279\times 10^{-11}$ & $6.352\times 10^{-14}$ \\
LGVI& $1.0$ & $40000$ & $335$ & - & $1.512\times 10^{-8}$ & $3.991\times 10^{-12}$ & $4.786\times 10^{-14}$ \\
    \hline\hline
   \end{tabular}
   \begin{tablenotes}
    \item[*] $h$ is integration step size, in seconds, fixed for LGVI but averaged over the run's duration for RKF7(8)
    \item[\ensuremath{\star}] $N_u$ is the total number of calculations of the mutual potential derivatives made within the run
    \item[\ensuremath{\diamond}] $t_{W}$ is the "wall-clock" time to complete each simulation run, in seconds
    \item[\ensuremath{\triangleleft}] $\epsilon$ is the error tolerance for the variable step size in RKF7(8)
    \item[\ensuremath{\ddag}] $\mbox{TE}$ and $\pi_T$ are total energy and the total angular momentum, respectively, while $\Delta$ refers to deviation from the initial value over simulation
    \item[\ensuremath{\dagger}] $\mathrm{E}[\cdot]$ denotes mean
   \end{tablenotes}
  \end{threeparttable}
 \end{center}
\end{table}

Here we see that for any pair of simulations, one using the LGVI and the other using the RKF7(8) scheme, for which the total energy metric performs about the same, the computation time needed to complete the simulation using the LGVI is a fraction of that needed using the RKF7(8). Simultaneous with this improvement in run time, the total angular momentum and attitude error metrics still perform better in the LGVI run than in the RKF7(8) run by multiple orders of magnitude. Going in the other direction, as the step size for the LGVI is reduced so that the computational burden using it begins to approach that for any chosen run using the RKF7(8), all error metrics remain at the same level as or else orders of magnitude smaller than those for the chosen RKF7(8) run. For the LGVI, the round-off error accumulates when multiplying rotation matrices at \refeqn{updateR}. The rotation matrix error of the LGVI is caused only by the floating-point arithmetic operation, and it is increased as the number of integration steps is increased. A similar trend is observed in the total angular momentum error for the LGVI, because determination of the total angular momentum in the inertial frame from the states written to the output file makes use of the rotation matrices.

\begin{figure}[ht]\centering
 \begin{subfigmatrix}{2}
  \subfigure[Trajectory of binary octahedra system components]{\includegraphics[width=0.50\textwidth]{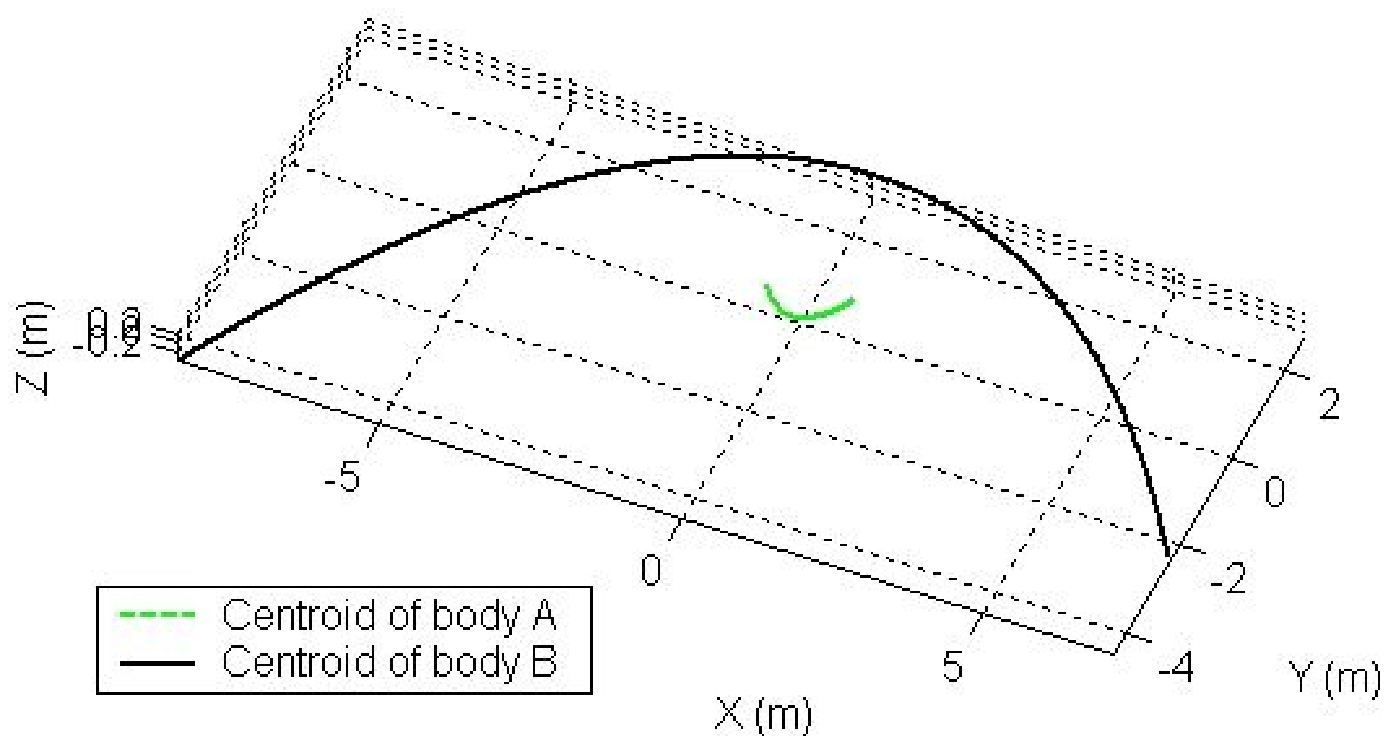}\label{fig:scenario3path}}
  \subfigure[Eccentricity and semi major axis]{\includegraphics[width=0.40\textwidth]{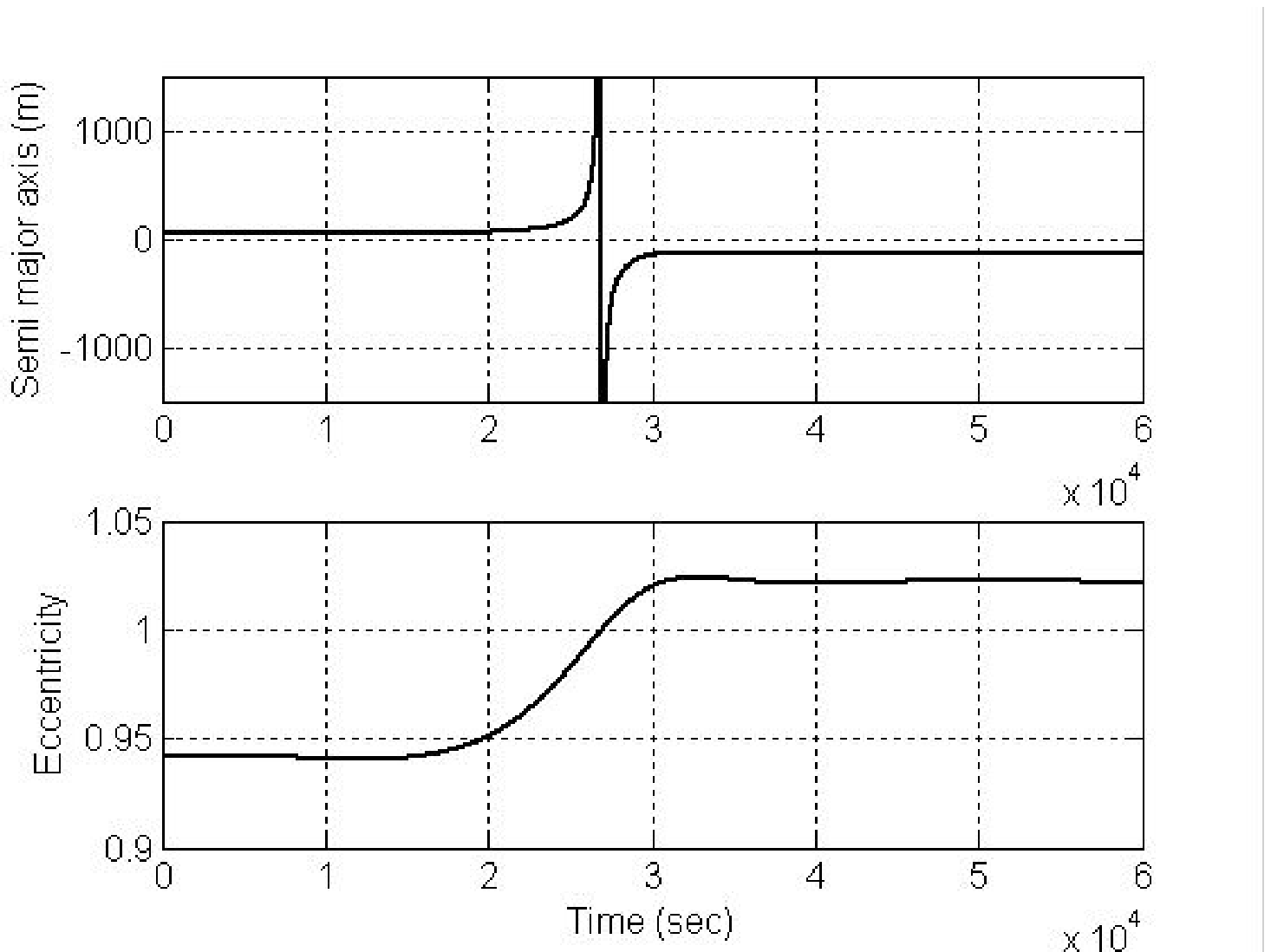}\label{fig:scenario3eccanda}}
\end{subfigmatrix}
\caption{(Scenario 3) Disruption of the binary octahedra system.
\label{fig:scenario3}}
\end{figure}

\paragraph{Scenario 3}%
The next scenario illustrates the ability of our methods to capture
the interesting effects of coupling in a mutual orbit configuration
that the Keplerian two-body approximation incorrectly predicts as
being perpetual. Simulation with the LGVI yields the trajectory
illustrated in Figure~\reffig{scenario3path}, which transitions from
a highly elliptical orbit to a hyperbolic escape path. This is shown
by the plots in Figure~\reffig{scenario3eccanda} of the semi-major
axis and eccentricity change during the close encounter, which
occurs roughly midway through the run duration of 60,000 seconds.
The initial conditions and body configurations are symmetric about
the initial orbital plane, and as such the motion of the centroids
should be restricted to the initial orbital plane. This is observed
numerically, as the body centroids remain within $8.6\times
10^{-14}$ meters of the initial orbital plane throughout the
simulation.

\begin{figure}[t]\centering
\includegraphics[width=0.45\textwidth]{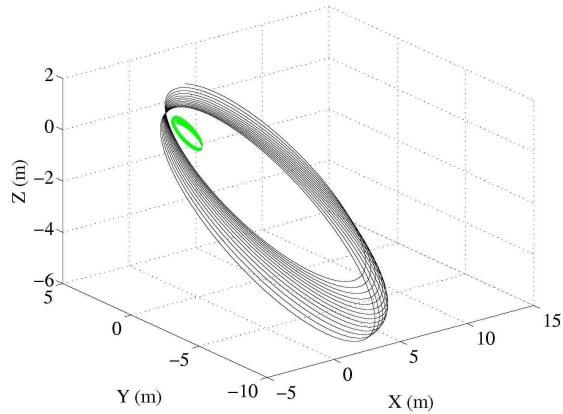}
\caption{(Scenario 4) Trajectory of binary octahedra
system}\label{fig:scenario4traj}
\end{figure}

\paragraph{Scenario 4} Finally, we examine a very long duration simulation starting from initial conditions that are stable in the sense that orbital trajectories are confined to separated and bounded regions but are highly unstable in the sense that body attitudes vary greatly and
irregularly. The trajectory of binary octahedra system is shown in
Figure~\reffig{scenario4traj}. Note that this motion cannot be
observed with the point mass assumption in the classical two body
problem.
This scenario is simulated both with LGVI at different step sizes
and with RKF7(8) at different error tolerances, for $5\times 10^6$
(sec) maneuver time. Table \reftab{scenario4perfcompare} illustrates
the general superiority of the LGVI approach, similar to Scenario 2.
In addition, we see that the simulation results of RKF7(8) with
larger error tolerances, $10^{-13}$ and $10^{-10}$, are completely
unreliable since the rotation matrix error is increased to the
unacceptable levels of $10^{-1}$ and $10^{1}$.

\begin{figure}[t]\centering
  \subfigure[Running time v.s. number mutual potential derivatives calculations]{\includegraphics[width=0.36\textwidth]{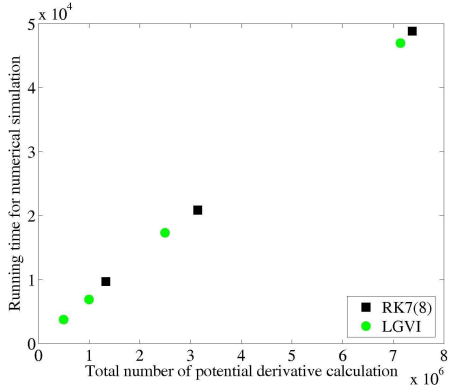}\label{fig:scenario4runtime}}\hspace*{2cm}
  \subfigure[Mean deviation of TE v.s. number of mutual potential derivatives calculations]{\includegraphics[width=0.39\textwidth]{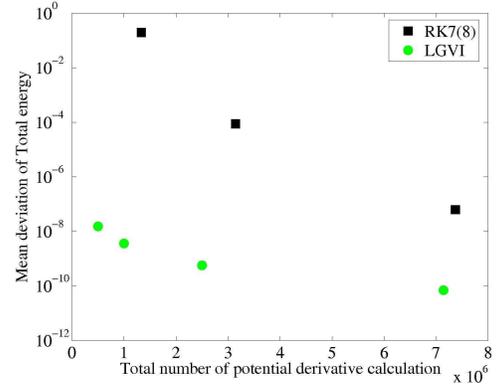}\label{fig:scenario4te}}
\caption{(Scenario 4) Quantitative comparisons of running time and
standard deviation of  total energy over number of mutual potential
derivatives calculations. \label{fig:scenario4quan}}
\end{figure}

\begin{figure}[h]\centering
  \subfigure[LGVI]{\includegraphics[width=0.39\textwidth]{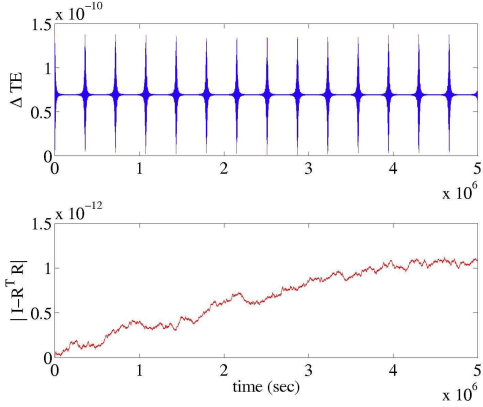}\label{fig:scenario4lgvi}}\hspace*{2cm}
  \subfigure[RKF7(8)]{\includegraphics[width=0.39\textwidth]{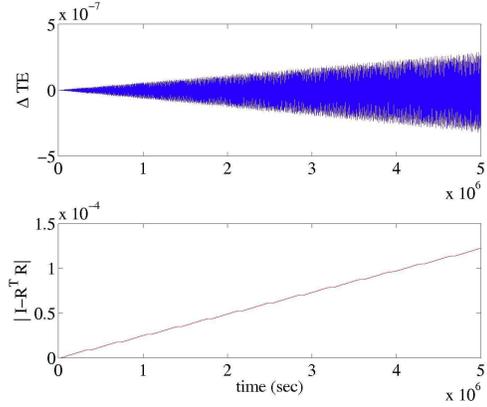}\label{fig:scenario4rk}}
\caption{(Scenario 4) Qualitative comparisons of behavior of
deviation in total energy and rotation matrix error over time with
similar computational load for LGVI and RKF7(8).
\label{fig:scenario4qual}}
\end{figure}

The quantitative comparison is summarized in Figure
\reffig{scenario4quan}, where the total computation time and the
standard deviation of the total energy are shown over the number of
the mutual potential derivatives calculations. In Figure
\reffig{scenario4runtime}, we see that the total computation time is
directly proportional to the number of the mutual potential
derivative calculations regardless of the integration scheme used.
This justifies the statement that the LGVI is an ``almost explicit''
computational method; the computational load to solve the implicit
equations is comparatively negligible. Considering that the mutual potential
derivative computations are the major computational burden for the
numerical simulation of full two rigid body dynamics, an efficient
integration scheme with fewer mutual potential derivative
calculations should exhibit smaller error measures. In Figure
\reffig{scenario4te}, the lower left portion represents higher
computational efficiency, and the upper right portion represents
lower efficiency. We see that LGVI is more efficient than RKF7(8).
Figure \reffig{scenario4qual} compares the total energy deviation
history and the relative rotation matrix error history. For LGVI,
the total energy varies within the level of $10^{-11}$. The repeated
peaks of the total energy deviation correspond to the perigees of
the orbit shown in Figure \reffig{scenario4traj}, but but there is
no tendency for the mean of the variation to drift for the entire
simulation time. The rotation matrix error is slightly increasing,
since the round off error in multiplying orthogornal matrices at
\refeqn{updateR} is accumulated. But the error is below $10^{-12}$
over the entire simulation's span. However, for RKF7(8), the
variation of the total energy is linearly increasing over time, and
the rotation matrix error is increasing to the level of $10^{-4}$.
We have seen that LGVI exhibits good total energy behavior for
exponentially long time periods, and it also preserves the group
structure well. The accuracy of RKF7(8) is vulnerable to a long time
simulation of the full two rigid body dynamics.

\begin{table}[t]
 \begin{center}
  \begin{threeparttable}
  \caption{(Scenario 4) Performance comparison between RKF7(8) and LGVI} \label{tab:scenario4perfcompare}
   \begin{tabular}{l|cccc|ccc}\hline\hline
Method & $h$ & $N_u$ & $t_{W}$ & $\epsilon$
& $\mathrm{E}[|\Delta\mbox{TE}|]$ & $\mathrm{E}[\norm{\Delta \pi_T}]$ & $\mathrm{E}[\norm{I-R^TR}]$ \\
\hline RKF7(8) & $8.81$ & $7373106$ & $48850$ & $10^{-16}$%
& $6.198\times10^{-8}$ & $3.830\times10^{-6}$ & $6.089\times 10^{-5}$\\%
RKF7(8) & $20.66$ & $3145532$ & $20873$ &  $10^{-13}$%
&$8.731\times10^{-5}$ & $7.427\times10^{-3}$ & $1.015\times 10^{-1}$\\%
RKF7(8) & $48.74$ & $1333475$ & $9699$ &  $10^{-10}$%
& $1.959\times10^{-1}$ & $5.246\times10^{0}$ & $1.834\times 10^{1}$\\%
\hline
LGVI& $0.7 $ & $7142858$ & $46963$ & -%
& $6.947\times 10^{-11}$ & $2.190\times 10^{-12}$ & $6.645\times 10^{-13}$\\%
LGVI& $2 $ & $2500000$ & $17320$ & -%
& $5.687\times 10^{-10}$ & $5.077\times 10^{-13}$ & $3.663\times 10^{-13}$\\%
LGVI& $5$ & $1000000$ & $6897$ & -%
& $3.601\times 10^{-9}$ & $8.447\times 10^{-13}$ & $1.642\times 10^{-13}$\\%
LGVI& $10$ & $500000$ & $3750$ & -%
& $1.517\times 10^{-8}$ & $2.567\times 10^{-13}$ & $3.521\times 10^{-13}$\\%
    \hline\hline
   \end{tabular}
   \begin{tablenotes}
    \item See Table~\reftab{perfcompare} for notations.
   \end{tablenotes}
  \end{threeparttable}
 \end{center}
\end{table}

\section{Conclusions}\label{sect:conclusions}

This paper presents a complementary combination of tools or
algorithms that achieves superior performance in the computationally
demanding simulation of the full two rigid bodies problem: a simple
method for efficient calculation of the mutual gravitational
potential and its derivatives given versatile polyhedral models of
the rigid bodies and a Lie group variational integrator consisting
of discrete relative equations of motion that preserves the
geometric features and structure of the configuration space. The use
of these easily implemented methods together allows for simulation
of the full two rigid bodies, capturing their fully coupled
dynamics, that is both accurate and efficient. The results above
show the accuracy maintained for energy, momenta, and attitude
geometry constraints, even for long simulation run times. The
presented mutual potential computations for polyhedral
approximations are themselves efficient, but for bodies of relevant
(i.e. large) model size, they still represent the bulk of the
computational burden during simulation, regardless of the
integration scheme used. Hence use of the Lie Group Variational
Integrator, which requires minimal mutual potential computations,
achieves an even greater combined efficiency. This yields maximal
simulation durations for fixed computational resources.

\section*{Acknowledgments}

The first author would like to acknowledge the support of the U. S.
Air Force Office of Scientific Research (AFOSR) for the period in
which this paper was written. The research of the third author was
partially supported by NSF grant DMS-0504747 and a University of
Michigan Rackham faculty grant.

\bibliography{aas}

\begin{thebibliography}{10}
\newcommand{\enquote}[1]{``#1''}

\bibitem{Bottke-Icarus-1996}
Bottke, W.~F. and Melosh, H.~J., \enquote{The Formation of Binary Asteroids and
  Doublet Craters,} {\em Icarus\/}, Vol.~124, 1996, pp.~372–391.

\bibitem{Bottke-Nature-1996}
Bottke, W.~F. and Melosh, H.~J., \enquote{The Formation of Asteroid Satellites
  and Doublet Craters by Planetary Tidal Forces,} {\em Nature\/}, Vol.~381,
  1996, pp.~51–53.

\bibitem{Margot-Science-2002}
Margot, J.~L., Nolan, M.~C., Benner, L. A.~M., Ostro, S.~J., Jurgens, R.~F.,
  Giorgini, J.~D., Slade, M.~A., and Campbell, D.~B., \enquote{Binary Asteroids
  in the {N}ear-{E}arth {O}bject Population,} {\em Science\/}, Vol.~296, May
  2002, pp.~1445--1448.

\bibitem{Merline-ASTIII-2002}
Merline, W.~J., Weidenschilling, S.~J., Durda, D.~D., Margot, J.~L., Pravec,
  P., and Storrs, A.~D., \enquote{Asteroids do have satellites,} {\em Asteroids
  {III}\/}, edited by W.~F. Bottke et~al., Space Science Series, University of
  Arizona, Tuscon, AZ, 2002, pp. 289--312.

\bibitem{Ostro-ASTIII-2002}
Ostro, S.~J., Hudson, R.~S., Benner, L. A.~M., Giorgini, J.~D., Magri, C.,
  Margot, J.-L., and Nolan, M.~C., \enquote{Asteroid Radar Astronomy,} {\em
  Asteroids {III}\/}, edited by W.~F. Bottke et~al., Space Science Series,
  University of Arizona, Tuscon, AZ, 2002, pp. 151--168.

\bibitem{Maciejewski-CMDA-1995}
Maciejewski, A.~J., \enquote{Reduction, Relative Equilibria and Potential in
  the Two Rigid Bodies Problem,} {\em Celestial Mechanics and Dynamical
  Astronomy\/}, Vol.~63, No.~1, 1995, pp.~1--28.

\bibitem{Scheeres-CMDA-2002}
Scheeres, D.~J., \enquote{Stability in the {F}ull {T}wo-{B}ody {P}roblem,} {\em
  Celestial Mechanics and Dynamical Astronomy\/}, Vol.~83, 2002, pp.~155--169.

\bibitem{Scheeres-NTAC-2003}
Scheeres, D.~J., \enquote{Stability of Relative equilibria in the {F}ull
  {T}wo-{B}ody {P}roblem,} {\em New Trends in Astrodynamics Conference\/}, Jan.
  2003.

\bibitem{Scheeres-AASastrospec-2003}
Scheeres, D.~J. and Augenstein, S., \enquote{Spacecraft motion about binary
  asteroids,} {\em Proc. {AAS/AIAA} Astrodynamics Specialist Conference\/},
  Aug. 2003.

\bibitem{Gabern-ICDSDE-2004}
Gabern, F., Koon, W.~S., and Marsden, J.~E., \enquote{Spacecraft dynamics near
  a binary asteroid,} {\em Proceedings of the fifth international conference on
  dynamical systems and differential equations\/}, Jun 2004.

\bibitem{Scheeres-DS-2005}
Scheeres, D.~J. and Bellerose, J., \enquote{The {R}estricted {H}ill {F}ull
  4-{B}ody {P}roblem: application to spacecraft motion about binary asteroids,}
  {\em Dynamical Systems: An International Journal\/}, Vol.~20, No.~1, 2005,
  pp.~23--44.

\bibitem{Borderies-celesmech-1978}
Borderies, N., \enquote{Mutual gravitational potential of N solid bodies,} {\em
  {C}elestial Mechanics\/}, Vol.~18, No.~3, 1978, pp.~295--307.

\bibitem{Braun-thesis-1991}
Braun, C.~V., {\em The gravitational potential of two arbitrary, rotating
  bodies with applications to the Earth-Moon system\/}, Ph.D. thesis,
  University of Texas at Austin, 1991.

\bibitem{Moritz-book-1980}
Moritz, H., {\em Advanced Physical Geodesy\/}, Abacus Press, 1980.

\bibitem{Geissler-Icarus-1996}
Geissler, P., Petit, J.-M., Durda, D.~D., Greenberg, R., Bottke, W., Nolan, M.,
  and Moore, J., \enquote{Erosion and Ejecta Reaccretion of 243 {I}da and Its
  Moon,} {\em Icarus\/}, Vol.~120, No.~1, 1996, pp.~140--157.

\bibitem{Ashenberg-JGCD-2005}
Ashenberg, J., \enquote{Proposed Method for Modeling the Gravitational
  Interaction Between Finite Bodies,} {\em Journal of Guidance, Control, and
  Dynamics\/}, Vol.~28, No.~4, 2005, pp.~768--774.

\bibitem{Werner-CMDA-1997}
Werner, R.~A. and Scheeres, D.~J., \enquote{Exterior Gravitation of a
  Polyhedron Derived and Compared with Harmonic and Mascon Gravitation
  Representations of Asteroid 4769 {C}astalia,} {\em Celestial Mechanics and
  Dynamical Astronomy\/}, Vol.~65, No.~3, 1997, pp.~313--344.

\bibitem{Werner-CMDA-2005}
Werner, R.~A. and Scheeres, D.~J., \enquote{Mutual Potential of Homogenous
  Polyhedra,} {\em Celestial Mechanics and Dynamical Astronomy\/}, Vol.~91,
  No.~3, March 2005, pp.~337--349.

\bibitem{Fahnestock-AASastrospec-2005}
Fahnestock, E.~G., Scheeres, D.~J., McClamroch, N.~H., and Werner, R.~A.,
  \enquote{Simulation and Analysis of Binary Asteroid Dynamics Using Mutual
  Potential and Potential Derivatives Formulation,} {\em Proc. {AAS/AIAA}
  Astrodynamics Specialist Conference\/}, Lake Tahoe, CA, Aug. 2005.

\bibitem{Fahnestock-CMDA-2006}
Fahnestock, E.~G. and Scheeres, D.~J., \enquote{Simulation of the Full Two
  Rigid Body Problem Using Polyhedral Mutual Potential and Potential
  Derivatives Approach,} {\em Celestial Mechanics and Dynamical Astronomy\/},
  submitted for publication.

\bibitem{Hairer-book-2000}
Hairer, E., Lubich, C., and Wanner, G., {\em {G}eometric {N}umerical
  {I}ntegration\/}, Springer, 2000.

\bibitem{Marsden-actanum-2001}
Marsden, J.~E. and West, M., \enquote{Discrete mechanics and variational
  integrators,} {\em Acta Numerica\/}, Vol.~10, 2001, pp.~357--514.

\bibitem{Iserles-actanum-2000}
Iserles, A., Munthe-Kaas, H.~Z., N{\o}rsett, S.~P., and Zanna, A.,
  \enquote{Lie-group methods,} {\em Acta Numerica\/}, Vol.~9, 2000,
  pp.~215--365.

\bibitem{Lee-IEEECCA-2005}
Lee, T., Leok, M., and McClamroch, N.~H., \enquote{A {L}ie group variational
  integrator for the attitude dynamics of a rigid body with application to the
  3{D} pendulum,} {\em Proceedings of the {IEEE} {C}onference on {C}ontrol
  {A}pplication\/}, Toronto, Canada, Aug. 2005, pp. 962--967.

\bibitem{Lee-CMAME-2005}
Lee, T., Leok, M., and McClamroch, N.~H., \enquote{{L}ie group variational
  integrators for the {F}ull {B}ody {P}roblem,} {\em {C}omputer {M}ethods in
  {A}pplied {M}echanics and {E}ngineering\/}, submitted, {A}vailable:
  \url{http://arxiv.org/abs/math.NA/0508365}.

\bibitem{Yeomans-Science-2000}
Yeomans, D.~K. et~al., \enquote{Radio Science Results During the
  {NEAR}-{S}hoemaker Spacecraft Rendezvous with {E}ros,} {\em Science\/},
  Vol.~289, Sept. 2000, pp.~2085--2088.

\bibitem{Fujiwara-Science-2006}
Fujiwara, A., Kawaguchi, J., Yeomans, D.~K., Abe, M., Mukai, T., Okada, T.,
  Saito, J., Yano, H., Yoshikawa, M., Scheeres, D.~J., Barnouin-Jha, O., Cheng,
  A.~F., Demura, H., Gaskell, R.~W., Hirata, N., Ikeda, H., Kominato, T.,
  Miyamoto, H., Nakamura, A.~M., Nakamura, R., Sasaki, S., , and Uesugi, K.,
  \enquote{The Rubble-Pile Asteroid Itokawa as Observed by Hayabusa,} {\em
  Science\/}, accepted for publication.

\bibitem{Richardson-Icarus-1998}
Richardson, D.~C., William F.~Bottke, J., and Love, S.~G., \enquote{Tidal
  Distortion and Disruption of Earth-Crossing Asteroids,} {\em Icarus\/},
  Vol.~134, 1998, pp.~47--76.

\bibitem{Roig-Icarus-2003}
Roig, F., Duffard, R., Penteado, P., Lazzaro, D., and Kodama, T.,
  \enquote{Interacting ellipsoids: a minimal model for the dynamics of
  rubble-pile bodies,} {\em Icarus\/}, Vol.~165, No.~2, 2003, pp.~355--370.

\bibitem{Korycansky-ASS-2004}
Korycansky, D.~G., \enquote{Orbital Dynamics for Rigid Bodies,} {\em
  Astrophysics and Space Science\/}, Vol.~291, 2004, pp.~57--74.

\bibitem{Duan-JCC-2001}
Duan, Z.-H. and Krasny, R., \enquote{An Adaptive Treecode for Computing
  Nonbonded Potential Energy in Classical Molecular Systems,} {\em Journal of
  Computational Chemistry\/}, Vol.~22, No.~2, 2001, pp.~184--195.

\end{thebibliography}
\bibliographystyle{aiaa}

\end{document}